\documentclass{amsart}
\usepackage{color}

\usepackage{amsmath} 
\usepackage{amssymb}
\usepackage{mathrsfs}

\newtheorem{theorem}{Theorem}[section] 
\newtheorem{claim}[theorem]{Claim}

\newtheorem{conclusion}[theorem]{Conclusion}
\newtheorem{observation}[theorem]{Observation}

\theoremstyle{definition}
\newtheorem{definition}[theorem]{Definition}

\newtheorem{fact}[theorem]{Fact}

\theoremstyle{remark}
\newtheorem{remark}[theorem]{Remark}
\newtheorem{question}[theorem]{Question}
\newtheorem{notation}[theorem]{Notation}

\newcommand{\Th}{{\rm Th}}

\newcommand{\trp}{{\rm trp}}

\newcommand{\tr}{{\rm tr}}
\newcommand{\fr}{{\rm fr}}

\newcommand{\cc}{{\rm cc}}
\newcommand{\wa}{{\rm wa}}
\newcommand{\car}{{\rm car}}

\newcommand{\lcf}{{\rm lcf}}
\newcommand{\TV}{{\rm TV}}

\newcommand{\ZFC}{{\rm ZFC}}
\newcommand{\GCH}{{\rm GCH}}
\newcommand{\Hom}{{\rm Hom}}

\newcommand{\com}{{\rm com}}
\newcommand{\comp}{{\rm comp}}

\newcommand{\hor}{{\rm hor}}
\newcommand{\ver}{{\rm ver}}
\newcommand{\spec}{{\rm spec}}
\newcommand{\upf}{{\rm upf}}

\newcommand{\uf}{{\rm uf}}
\newcommand{\Card}{{\rm Card}}

\newcommand{\at}{{\rm at}}
\newcommand{\Ord}{{\rm Ord}}

\newcommand{\id}{{\rm id}}

\newcommand{\Dom}{{\rm Dom}}
\newcommand{\Rang}{{\rm Rang}}

\newcommand{\bfm}{{\mathbf m}}

\newcommand{\bfh}{{\mathbf h}}

\newcommand{\bfk}{{\mathbf k}}

\newcommand{\rest}{{\restriction}}

\newcommand{\wilog}{{\rm without loss of generality}}
\newcommand{\Wilog}{{\rm Without loss of generality}}

\newcommand{\then}{{\underline{then}}}
\newcommand{\when}{{\underline{when}}}
\newcommand{\where}{{\underline{where}}}

\newcommand{\Iff}{{\underline{iff}}}
\newcommand{\mn}{{\medskip\noindent}}
\newcommand{\sn}{{\smallskip\noindent}}

\newcommand{\bbB}{{\mathbb B}}

\newcommand{\cC}{{\mathscr C}}

\newcommand{\cE}{{\mathscr E}}

\newcommand{\cF}{{\mathscr F}}

\newcommand{\cI}{{\mathscr I}}

\newcommand{\bbL}{{\mathbb L}}
\newcommand{\bbN}{{\mathbb N}}

\newcommand{\cP}{{\mathscr P}}

\newcommand{\cT}{{\mathscr T}}
 
\newcommand{\varp}{{\varepsilon}}

 \newcommand{\cf}{{\rm cf}}

\newcount\skewfactor
\def\mathunderaccent#1#2 {\let\theaccent#1\skewfactor#2
\mathpalette\putaccentunder}
\def\putaccentunder#1#2{\oalign{$#1#2$\crcr\hidewidth
\vbox to.2ex{\hbox{$#1\skew\skewfactor\theaccent{}$}\vss}\hidewidth}}

\newenvironment{PROOF}[2][\proofname.]
   {\begin{proof}[#1]\renewcommand{\qedsymbol}{\textsquare$_{\rm #2}$}}
   {\end{proof}}

\usepackage[hidelinks]{hyperref}

\begin{document}
\makeatletter\def\shfiuwefootnote{\gdef\@thefnmark{}\@footnotetext}\makeatother\shfiuwefootnote{Version 2018-08-06\_12. See \url{https://shelah.logic.at/papers/1026/} for possible updates.}

\title[The spectrum of ultraproducts]{The spectrum of 
ultraproducts of finite \\
cardinals for an
  ultrafilter \\ 
Sh1026}
\author {Saharon Shelah}
\address{Einstein Institute of Mathematics\\
Edmond J. Safra Campus, Givat Ram\\
The Hebrew University of Jerusalem\\
Jerusalem, 91904, Israel\\
 and \\
 Department of Mathematics\\
 Hill Center - Busch Campus \\ 
 Rutgers, The State University of New Jersey \\
 110 Frelinghuysen Road \\
 Piscataway, NJ 08854-8019 USA}
\email{shelah@math.huji.ac.il}
\urladdr{http://shelah.logic.at}
\thanks{The author would like to thank the ISF for partial support of
  this research, Grant No. 1053/11.  The reader should note that the
  version in my website is usually more updated than the one in the
  mathematical arXive.  
The author thanks Alice Leonhardt for the beautiful typing.
  First typed April 5, 2013.}





\date{August 6, 2018}

\begin{abstract}
We complete the characterization of the possible spectrum of regular
ultrafilters $D$ on a set $I$, where the spectrum is the set of 
ultraproducts of (finite) cardinals modulo $D$ which are infinite.
\end{abstract}

\maketitle
\numberwithin{equation}{section}
\setcounter{section}{-1}
\newpage

\section {Introduction}
\bigskip

\subsection {Background, questions and results}\
\bigskip

Ultraproducts were very central in model theory in the sixties, usually for
regular ultrafilters.  The question of ultraproducts of infinite
cardinals had been resolved (see \cite{CK73}): letting $D$ be a 
regular ultrafilter on a set $I$, (for transparency we ignore the case
of a filter)
\mn
\begin{enumerate}
\item[$(*)_1$]  if $\bar\lambda = \langle \lambda_s:s \in
  I\rangle$ and $\lambda_s \ge \aleph_0$ for $s \in I$ \then \, 
$\prod\limits_{s \in I} \lambda_s/D = \mu^{|I|}$ when $\mu = 
\lim \sup_D(\bar\lambda) := \sup\{\chi$: the cardinal $\chi$ satisfies
$\{s \in I:\lambda_s \ge \chi\} \in D\}$.
\end{enumerate}
\mn
What about the ultraproducts of finite cardinals?  Of course, under
naive interpretation, if $\{\lambda_s:\lambda_s = 0\} \ne \emptyset$ 
the result is zero, so for notational
simplicity we always assume $s \in I \Rightarrow \lambda_s \ge 1$.
Also for every $n \ge 1$, letting $\lambda_s = n$ for $s \in I$ we have
$\prod\limits_{s} \lambda_s/D=n$ so the question was
\begin{question}  
\label{y4}
Given an infinite set $I$
\mn
\begin{enumerate}
\item[$(a)$]  [the singleton problem] what infinite cardinals $\mu$
belong to $\cC_I = \cC^{\car}_I$, i.e. can be represented as 
$\{\prod\limits_{s \in I} \lambda_s/D:D$ a regular ultrafilter on
  $I,1 \le \lambda_s < \aleph_0\} \backslash \{\lambda:1 \le \lambda <
  \aleph_0\}$
\sn
\item[$(b)$]  [the spectrum problem] moreover what are the possible 
spectra, i.e. which sets of cardinals belong to $\mathbf C_I$ which 
is the family of sets $\cC$ such that for some $D$, a regular ultrafilter on
$I$ we have $\cC = \upf(D)$ where $\upf(D) = 
\{\prod\limits_{s \in I} \lambda_s/D:1 \le \lambda_s < \aleph_0$ for
$s \in I\} \backslash \{\lambda:1 \le \lambda < \aleph_0\}$
\end{enumerate}
\mn
Keisler \cite{Ke67a} asks and has started on \ref{y4}: 
(assuming GCH was prevalent at the time as the
situation was opaque otherwise)
\mn
\begin{enumerate}
\item[$(*)_2$]  assume GCH, a sufficient condition for $\cC \in \mathbf C_I$ is:
\sn
\begin{enumerate}
\item[$(a)$]  $\cC$ is a set of successor (infinite) cardinals
\sn
\item[$(b)$]  $\max(\cC) = |I|^+$
\sn
\item[$(c)$]  if $\mu = \sup\{\chi < \mu:\chi \in \cC\}$ then $\mu^+
  \in \cC$
\sn
\item[$(d)$]  if $\mu^+ \in \cC$ then $\mu \cap \cC$ has cardinality $<
 \mu$.
\end{enumerate}
\end{enumerate}
\mn
Keisler used products and $D$-sums of ultrafilters.  Concerning the 
problem for singletons a conjecture of Keisler \cite[bottom of
pg.49]{Ke67a} was resolved in \cite{Sh:7}:
\mn
\begin{enumerate}
\item[$(*)_3$]  $\mu = \mu^{\aleph_0}$ \when \, $\mu \in \cC_I$, i.e.
when $\mu = \prod\limits_{s \in I} \lambda_s/D$ is infinite, $D$ an
 ultrafilter on $I$, each $\lambda_s$ finite non-zero
\end{enumerate}
\mn
The proof uses coding enough ``set theory" on the $n$'s and using the
model theory of the ultra-product.  This gives a necessary condition
(for the singleton version),
but is it sufficient?  This problem was settled in
\cite[Ch.V,\S3]{Sh:a} = \cite[Ch.VI,\S3]{Sh:c} proving that this is
also a sufficient condition 
(+ the obvious condition $\mu \le 2^{|I|}$), that is
\mn
\begin{enumerate}
\item[$(*)_4$]   $\mu \in \cC_I := 
\cup\{\cC:\cC \in \mathbf C_I\}$ \Iff \, $\mu = \mu^{\aleph_0} \le 2^{|I|}$.
\end{enumerate}
\mn
The constructions in \cite[Ch.VI,\S3]{Sh:a} = \cite[Ch.VI,\S3]{Sh:c},
use a family $\cF$ of functions with domain $I$ and a filter $D$ on $I$
such that $\cF$ is independent over $D$ (earlier Kunen used such
family $\cF \subseteq {}^\lambda \lambda$ for constructing a good
ultrafilter on $\lambda$ in ZFC; eliminating the use of an instance of
$\GCH$ in the proof of Keisler; earlier Engelking-Karlowicz proved the
existence of such $\cF$).  In particular in the construction in
\cite[Ch.VI,\S3]{Sh:a} of maximal such
filters and the Boolean Algebra $\bbB = \cP(\lambda)/D$ are 
central.  We decrease the family and increase $D$ during the
construction; specifically we
construct $\cF_\ell(\ell \le n)$ decreasing with $\ell,D_\ell$ a 
filter on $I$ increasing with $\ell,D_\ell$ a maximal filter such that
$\cF_\ell$ is independent $\mod D_\ell$; so if $\cF_n = \emptyset$
then $D_0$ is an ultrafilter and we have 
$\bbB_\ell = \cP(I)/D_\ell$ is essentially $\lessdot$-decreasing and in the
ultrapowers $\bbN^I/D_\ell$ the part which $\bbB_\ell$ induces for $\ell
\le n$, is a sequence of initial
segments of $\bbN^\bbB/D_0$ decreasing with $\ell$.

In \cite[Ch.VI,Exercise 3.35]{Sh:a} = \cite[pg.370]{Sh:c} this is formalized:
\mn
\begin{enumerate}
\item[$(*)_5$]   if $D_0$ is a
filter on $I,\bbB_0 = \cP(I)/D_0,D_1 \supseteq D_0$ an ultrafilter, $D =
\{A/D_0:A \in D_1\}$ so $D \in \uf(\bbB_0)$ then $\bbN^{\bbB_0}/D^+_0$
is an initial segment of $\bbN^I/D$; (also $\bbB$ satisfies the
c.c.c., but this is just to ensure $\bbB$ is complete, anyhow this
holds in all relevant cases here).
\end{enumerate}
\mn
It follows that we can replace $\cP(I)$ by a Boolean
Algebra $\bbB_1$ extending $\bbB_0$.  The Boolean Algebra related to
$\cF$ is the completion of the
Boolean Algebra generated by $\{x_{f,a}:f \in \cF,a \in
\Rang(f)\}$ freely except $x_{f,a} \cap x_{f,b} =0$ for $a \ne b \in
\Rang(f)$ and $f \in \cF$.  
So if $\Rang(f)$ is countable for every $f \in \cF$, the
Boolean Algebra satisfies the $\aleph_1$-c.c. (in fact, is free), this 
was used there to deal with $\lcf(\kappa,D)$ for $\kappa = \aleph_0$ (for $\kappa > \aleph_0$
we need $\Rang(f) = \kappa$) and is continued lately 
in works of Malliaris-Shelah.
But for $\upf(D)$ only the case of $f$'s with countable range is
used.  

The problem of the spectrum (i.e. \ref{y4}(b)) was not needed in
\cite[Ch.VI,\S3]{Sh:a} for the
model theoretic problems which were the aim of \cite[Ch.VI]{Sh:a}, still
the case of finite spectrum was resolved there (also cofinality, i.e. 
$\lcf(\kappa,D)$ was addressed).

This was continued by Koppelberg \cite{Kp80} using a
possibly infinite $\lessdot$-increasing chains of complete Boolean
Algebras; also she uses a system of projections instead of maximal filters but
this is a reformulation as this is equivalent, see \ref{x23} below.

Koppelberg \cite{Kp80} returns to the full spectrum problem proving:
\mn
\begin{enumerate}
\item[$(*)_6$]   $\cC \in \mathbf C_I$ \when \, $\cC$ satisfies:
\sn
\begin{enumerate}
\item[$(a)$]  $\cC \subseteq \Card$
\sn
\item[$(b)$]  $\max(\cC) = 2^{|I|}$
\sn
\item[$(c)$]  $\mu = \mu^{\aleph_0}$ if $\mu \in \cC$
\sn
\item[$(d)$]  if $\mu_n \in \cC$ for $n < \omega$ then
  $\prod\limits_{n} \mu_n \in \cC$.
\end{enumerate}
\end{enumerate}
\mn
Central in the proof is $(*)_5$ above
(\cite[Ch.VI,Ex3.35,pg.370]{Sh:a}).  The result of Koppelberg 
is very strong, still the full characterization is not 
obtained; also Kanamori in his math review of her work asked about it.
\end{question}

\noindent
Here we give a complete answer to the spectrum problem \ref{y4}(b),
that is, Theorem \ref{k6} gives a full $\ZFC$ answer to \ref{y4}, that
is.
\begin{theorem}
\label{y8}
For any infinite set $I,\cC \in \mathbf C_I$ \Iff \, $\cC$ is a set of
cardinals such that $\mu \in \cC \Rightarrow \mu = \mu^{\aleph_0} \le
2^{|I|}$ and $2^{|I|} \in \cC$.
\end{theorem}

\noindent
We now comment on some further questions on ultra-powers.  

The problem of cofinalities was central in \cite[Ch.VI,\S3]{Sh:a} in
particular, $\lcf(\aleph_0,\lambda))$ (see \ref{x16} below).
[Why?  E.g. if $\Th(M)$, the complete first order theory of the model $M$
 is unstable then $M^I/D$ is not $\lcf(\aleph_0,\lambda)^+$-saturated.]  
Another question was raised by the author \cite[pg.97]{Sh:14} and
 independently by Eklof \cite{Ek73}:

\begin{question}  
\label{y9}
Assume $f_n \in {}^I\bbN,f_{n+1} <_D f_n$ and $\mu \le \prod\limits_{s
  \in I} f_n(s)/D$ for every $n$ \then \, is there
$f \in {}^I\bbN$ such that $f <_D f_n$ for every 
$n$ and $\mu \le \prod\limits_{s \in I} f(s)/D$?

The point in \cite[pg.75]{Sh:14} was investigating saturation of
ultrapowers (and ultraproducts) and Keisler order on first order
theories.  The point in \cite{Ek73} was ultraproduct of Abelian groups.

To explain the cofinalities problem, see \ref{y12}.   
We can consider the following: 
for $D$ a regular ultrafilter on $I$ we consider $M =
\bbN^\lambda/D$; for $a \in M$ let $\lambda_a = |\{b:b <_M a\}|$ and
define $E_M = \{(a,b):a,b \in M$ and $\lambda_a = \lambda_b \ge
\aleph_0\}$.  So $E_M$ is a convex equivalence relation, and the equivalence
classes are naturally linearly ordered and let $A_{D,\lambda} = \{a \in
M:\lambda_a = \lambda\}$.  So $\upf(D) = \{\lambda_a:A_{D,\lambda} \ne
\emptyset\}$ and Question \ref{y9} asks: can the co-initiality of some
$A_{D,\lambda}$ be $\aleph_0$.  As $M$ is $\aleph_1$-saturated, in
this case the cofinality of $M \rest \{c:\lambda_c < \lambda_a$ (hence
$c <_M a)\}$ is $\lcf(\aleph_0,D)$ which is the co-initiality of
$A_{D,\min(\upf(D))}$. 
\end{question}

\noindent
So a natural question is
\begin{question}
\label{y12}
What are the possible $\spec_1(D) = \{(\lambda,\theta,\partial):\lambda \in
\upf(D),\partial$ the cofinality of $A_{D,\lambda}$ and $\theta$ the
co-initiality of $A_{D,\lambda}\}$ for $D$ a regular ultrafilter on $I$?  
\end{question}

\noindent
A further question is:
\begin{question}
\label{y15}
Assume $\kappa = \cf(\kappa) < \lambda_1 = \lambda^{\aleph_0}_1 < \lambda_2 =
\lambda^{\aleph_0}_2,\lambda^{<\kappa>_{\tr}}_1 \le 2^\lambda$; see
\ref{x18}(4).  Is
there a regular ultrafilter $D$ on $\lambda$ such that for $n_i \in \bbN$
for $i < \lambda$ we have $\prod\limits_{i} n_i/D = \lambda_1$ and 
$\prod\limits_{i} 2^{n_i}/D = \lambda_2$?

This work was presented in the May 2013 Eilat Conference honoring Mati
Rubin's retirement.  In a work in preparation \cite{Sh:F1427}, we try
to build a counterexample to question \ref{y9}.
\end{question}
\bigskip

\subsection {Preliminaries}\
\bigskip

We define $\lcf(\kappa,D)$ and $M^{\bbB}/D$, when $\bbB$ is a Boolean
Algebra and more.

\begin{definition}  
\label{x16}
For $D$ an ultrafilter on $I,\kappa$ a regular cardinal let $\mu =
\lcf(\kappa,D)$ be the co-initiality of the linear order $(\kappa^I/D)
\rest \{f/D:f \in {}^I \kappa$ is not $D$-bounded by any $\varepsilon <
\kappa\}$. 
\end{definition}

\begin{notation}
\label{x18}
1) $\bbB$ denotes a Boolean Algebra, usually complete; let
   $\comp(\bbB)$ be the completion of $\bbB$.

\noindent
2) $\uf(\bbB)$ is the set of ultrafilters on $\bbB$.

\noindent
3) Let $\bbB^+ = \bbB \backslash \{0_{\bbB}\}$. 

\noindent
4) Let $\cc(\bbB) = \min\{\kappa:\bbB$ satisfies the
   $\kappa$-c.c.$\}$, necessarily a regular cardinal.

\noindent
5) For $\lambda \ge \kappa = \cf(\kappa)$ let $\lambda^{\langle
  \kappa\rangle} = \trp_\kappa(\lambda) = \sup\{|\lim_\kappa(\cT)|:\cT
\subseteq {}^{\kappa >}\lambda$ is a subtree of cardinality $\le
\lambda\}$ where $\lim_\kappa(\cT) = \{\eta \in {}^\kappa \Ord:\eta
\rest i \in \cT$ for every $i < \kappa\}$.

\noindent
6) For a Boolean algebra $\bbB$ let $\comp(\bbB)$ be its completion.
\end{notation}

\begin{definition}  
\label{x20}
For a Boolean Algebra $\bbB$ and a model or a set $M$.

\noindent
1) Let $M^{\bbB}$ be the set of partial functions $f$ from
   $\bbB^+$ into $M$ such that for some maximal antichain $\langle
   a_i:i < i(*)\rangle$ of $\bbB,\Dom(f)$ includes $\{a_i:i < i(*)\}$
 and is included in\footnote{for the $D_\ell \in \uf(\bbB_\ell)$ 
ultra-product, \wilog \, $\bbB$ is complete, then
\wilog \, $f \rest \{a_i:i < i(*)\}$ is one to one.} $\{a \in \bbB^+
:(\exists i)(a \le a_i)\}$ and $f \rest \{a \in \Dom(f):
a \le a_i\}$ is constant for each $i$.

\noindent
1A) Naturally for $f_1,f_2 \in M^{\bbB}$ we say $f_1,f_2$ are
$D$-equivalent, or $f_1 = f_2 \mod D$ \when \, for some $b \in D$ we have
$a_1 \in \Dom(f_1) \wedge a_2 \in \Dom(f_2) \wedge a_1 \cap a_2 \cap b
 > 0_{\bbB} \Rightarrow f_1(a_1) = f_2(a_2)$.

\noindent
1B) Abusing notation, not only $M^{\bbB_1} \subseteq M^{\bbB_2}$ but
$M^{\bbB_1}/D_1 \subseteq M^{\bbB_2}/D_2$ when $\bbB_1 \lessdot
\bbB_2,D_\ell \in \uf(\bbB_\ell)$ for $\ell=1,2$ and $D_1 \subseteq
D_2$, that is, for $f \in M^{\bbB_1}$ we identify $f/D_1$ and $f/D_2$.

\noindent
2) For $D$ an ultrafilter on the completion of the Boolean Algebra
   $\bbB$ we define $M^{\bbB}/D$ naturally,
   as well as $\TV(\varphi(f_0,\dotsc,f_{n-1})) \in \comp(\bbB)$ \when \,
   $\varphi(x_0,\dotsc,x_{n-1}) \in \bbL(\tau_M)$ and $f_0,\dotsc,f_{n-1} \in
   M^{\bbB}$ where $\TV$ stands for truth value and $M^{\bbB}/D \models
\varphi[f_0/D,\dotsc,f_{n-1}/D]$ iff 
$\TV_M(\varphi(f_0,\dotsc,f_{n-1})) \in D$. 

\noindent
3) We say $\langle a_n:n < \omega\rangle \, D$-represents $f \in
\bbN^{\bbB}$ \when \, $\langle a_n:n < \omega\rangle$ is a maximal
antichain of $\bbB$ (allowing $a_n = 0_{\bbB}$)
and for some $f' \in \bbN^{\bbB}$ which is 
$D$-equivalent to $f$ (see \ref{x20}(1A)) we have $f'(a_n)=n$.  We
may omit $D$ if $D = \{1_{\bbB}\}$ and say just $\langle a_n:n <
\omega\rangle$ represents $f$.

\noindent
4) We say $\langle (a_n,k_n):n < \omega\rangle$ represents $f \in
   \bbN^{\bbB}$ \when \,:
\mn
\begin{enumerate}
\item[$(a)$]  the $k_n$ are natural numbers with no repetition
\sn
\item[$(b)$]  $\langle a_n:n < \omega \rangle$ is a maximal antichain
  of $\bbB$
\sn
\item[$(c)$]  $f(a_n) = k_n$.
\end{enumerate}
\mn
The proofs in \cite[Ch.VI,\S3]{Sh:a} use downward induction on the
cardinals.  
\end{definition}

\begin{observation}  
\label{x21}
1) If $\bbB$ is a complete Boolean Algebra and $f \in \bbN^{\bbB}$ \then
\, some sequence $\langle a_n:n < \omega\rangle$ represents $f$.  

\noindent
1A) If $\bbB$ is a c.c.c. Boolean Algebra and $f \in \bbN^{\bbB}$,
\then \, some sequence $\langle a_n,k_n:n < \omega\rangle$ represents $f$.

\noindent
2) For a model $M$ and Boolean Algebra $\bbB_1$ and ultrafilter $D$ on
 its completion $\bbB_2$ we have $M^{\bbB_1}/D = M^{\bbB_2}/D$.
\end{observation}

\begin{fact}  
\label{x30}
1) If $\bbB_1 \subseteq \bbB_2$ are Boolean Algebras, $\bbB$ is a
   complete Boolean Algebra and $\pi_1$ is a homomorphism from
   $\bbB_1$ into $\bbB$ \then \, there is a homomorphism $\pi_2$
   from $\bbB_2$ into $\bbB$ extending $\pi_1$.

\noindent
2) There is a homomorphism $\pi_3$ from $\bbB_3$ into $\bbB$ extending
$\pi_\ell$ for $\ell=0,1,2$ \when \,:
\mn
\begin{enumerate}
\item[$(a)$]  $\bbB_0 \subseteq \bbB_\iota \subseteq \bbB_3$ are
  Boolean Algebras for $\iota = 1,2$
\sn
\item[$(b)$]  $\bbB_1,\bbB_2$ are freely amalgamated over $\bbB_0$
  inside $\bbB_3$
\sn
\item[$(c)$]  $\bbB$ is a complete Boolean Algebra
\sn
\item[$(d)$]  $\pi_\ell$ is a homomorphism from $\bbB_\ell$ into
  $\bbB$ for $\ell=0,1,2$
\sn
\item[$(e)$]  $\pi_0 \subseteq \pi_1$ and $\pi_0 \subseteq \pi_2$.
\end{enumerate}
\end{fact}

\begin{PROOF}{\ref{x30}}  
1) Well known.

\noindent
2) Straightforward.
\end{PROOF}

\begin{observation}  
\label{x23}
Assume $\bbB_1 \lessdot \bbB_2$ are Boolean Algebras and $\bbB_1$ is complete.

\noindent
1) The following properties of $D$ are equivalent:
\mn
\begin{enumerate}
\item[$(a)$]  $D$ is a maximal filter on $\bbB_2$ (among those)
  disjoint to $\bbB_1 \backslash \{1_{\bbB_1}\}$
\sn
\item[$(b)$]  there is a projection $\pi$ of $\bbB_2$ onto $\bbB_1$
  such that $D = \{a \in \bbB_2:\pi(a) = 1_{\bbB_1}\}$.
\end{enumerate}
\mn
1A) Moreover $D$ determines $\pi$ uniquely and vice versa, in
particular $\pi(c)$ is the unique $c' \in \bbB_1$ such that $c=c' \mod D$.

\noindent
2) If $D$ satisfies (1)(a) and $D_1$ is an ultrafilter of $\bbB_1$,
\then \, there is a one and only one ultrafilter $D_2 \in \uf(\bbB_2)$
extending $D_1 \cup D$.
\end{observation}

\begin{PROOF}{\ref{x23}}
1) \underline{Clause (a) implies clause (b)}:

As $D$ is a filter on $\bbB_2$ clearly for some Boolean Algebra $\bbB'_2$, 
there is a homomorphism $\mathbf j_0:\bbB_2 \rightarrow \bbB'_2$ which
is onto, such that $a \in \bbB_2 \Rightarrow (a \in D \leftrightarrow
\mathbf j_0(a) = 1_{\bbB'_2})$.  As $D \cap \bbB_1 = \{1_{\bbB_1}\}$
necessarily $\mathbf j_0 \rest \bbB_1$ is one-to-one.  Let $\bbB'_1 =
\mathbf j_0(\bbB_1)$ so $\mathbf j_1 := (\mathbf j_0 \rest \bbB_1)^{-1}$ is an
isomorphism from $\bbB'_1$ onto $\bbB_1$ hence by \ref{x30}(1) and the
assumption that $\bbB_1$ is complete there
is a homomorphism $\mathbf j_2$ from $\bbB'_2$ onto $\bbB_1$ extending
$\mathbf j_1$.  Hence $\mathbf j_3 = \mathbf j_2 \circ \mathbf j_0$ is a
homomorphism from $\bbB_2$ onto $\bbB_1$ extending $\id_{\bbB_1}$, so
it is a projection.

Lastly, $\mathbf j^{-1}_3\{1_{\bbB_1}\}$ is a filter extending $D$ and
disjoint to $\bbB_1 \backslash \{1_{\bbB_1}\}$.  By the maximality of
$D$ we have equality.

An alternative proof is:

Let $\bbB'_2$ be the sub-algebra of $\bbB_2$ generated by $\bbB_1 \cup
D$.  Clearly every member of $\bbB'_2$ can be represented as $(a \cap
b) \cup ((1-a) \cap \sigma(\bar a,\bar b))$ with $a,a_m \in D$ for $m
< n = \ell g(\bar a)$ and $b \in \bbB_1,b_k \in \bbB_1$ for $k < \ell g(\bar
b),\sigma$ a Boolean term such that $\bigwedge\limits_{k<n} a \le a_k$,
equivalently $\bigwedge\limits_{k<n} a \cap (1-a_k) =0$.  
We try to define a function $\pi$ from $\bbB'_2$ into $\bbB_1$ by:
\mn
\begin{enumerate}
\item[$\oplus$]  $\pi((a \cap b) \cup ((1-a) \cap \sigma(\bar a,\bar
  b))) = b$ for $a,\bar a,b,\bar b$ as above.
\end{enumerate}
\mn
We have to prove that $\pi$ is as promised.
\mn
\begin{enumerate}
\item[$(*)_1$]  $\pi$ is a well defined (function from $\bbB'_2$ into
  $\bbB_1$). 
\end{enumerate}
\mn
Why?  Obviously for every $c \in \bbB'_2$ there are $a,\bar a,b,\bar
b,\sigma$ as above, so $\pi(c)$ has at least one definition, still we
have to prove that any two such definitions agree.  So 
assume $c = (a_\ell \cap b_\ell) \cup ((1-a_\ell) \cap
\sigma_\ell(\bar a_\ell,\bar b_\ell))$ for $\ell=1,2$ as above so
with $a_1,a_2,a_{1,k},a_{2,m} \in D$ and $b_1,b_2,\bar b_1,\bar b_2 
\in \bbB_1$ such that $a_1 \le a_{1,k}$ for every $k < \ell g(\bar
a_1),a_2 \le a_{2,m}$ for every $m < \ell g(\bar a_2)$.  
We should prove that $b_1 =
b_2$, if not \wilog \, $b_1 \nleq b_2$ hence $b := b_1 - b_2 >
0$.  Clearly $a := a_1 \cap a_2 \in D$ and computing $c \cap b \cap
a$ in two ways we get $a \cap b \cap b_1 = a \cap b \cap b_2$ hence $a \cap
b = a \cap b \cap b_1 = a \cap b \cap b_2 = a \cap 0=0$ recalling
$b=b_1-b_2$, hence $a
\le 1-b$ so as $a \in D$ necessarily $1-b \in D$.  But $b \in \bbB^+_1$
so $1-b \in \bbB_1 \backslash \{1_{\bbB_1}\}$, contradiction to the
assumption on $D$.
\mn
\begin{enumerate}
\item[$(*)_2$]   $\pi$ commutes with $``x \cap y"$.
\end{enumerate}
\mn
Why?  Assume that for $\ell=1,2$ we have $c_\ell = (a_\ell \cap
b_\ell) \cup ((1-a_\ell) \cap \sigma_\ell(\bar a_\ell,\bar b_\ell))$ 
with $a_\ell,b_\ell,\bar a_\ell,\bar b_\ell,\sigma_\ell$ as above.

So $\pi(c_\ell) = b_\ell$ and letting $a = a_1 \cap a_2 \in D$ we have
$c := c_1 \cap c_2 = (a \cap (b_1 \cap b_2)) \cup ((1-a) \cap
\sigma(\bar a,\bar b))$ where $\bar a = \bar a_1 \char 94 \langle
a_1 \rangle \char 94 \bar a_2 \char 94 \langle a_2 \rangle,\bar b =
\bar b_1 \char 94 \langle b_1 \rangle \char 94 \bar b_2 \char 94
\langle b_2\rangle$ for some suitable term $\sigma$.

As $a \in D$, clearly $\pi(c) = b_1 \cap b_2 = \pi(c_1) \cap
\pi(c_2)$, as required.
\mn
\begin{enumerate}
\item[$(*)_3$]   $\pi$ commutes with ``$1-x$".
\end{enumerate}
\mn
Why?  Let $c = (a \cap b) \cup ((1-a) \cap \sigma(\bar a,\bar b))$ 
hence $1-c= (a \cap (1-b)) \cup ((1-a) \cap (1- \sigma(\bar a,\bar
b))$ hence $\pi(1-c) = 1-b = 1 - \pi(c)$ so we are done.
\mn
\begin{enumerate}
\item[$(*)_4$]   $\pi$ is a projection onto $\bbB_1$.
\end{enumerate}
\mn
[Why?  By $(*)_1,(*)_2,(*)_3$ clearly $\pi$ is a homomorphism from
  $\bbB'_2$ into $\bbB_1$.  So its range is $\subseteq \bbB_1$ and 
if $c \in \bbB_1$ let $b=c,a = 1_{\bbB_1},\bar a = \langle \rangle 
= \bar b$ and $\sigma(\bar a,\bar b) = 0_{\bbB_1}$ so $c = (a \cap b)
\cup ((1-a) \cap \sigma(\bar a,\bar b)$ and 
$a,b,\bar a,\bar b,\sigma$ are as required so
  $\pi((a \cap b) \cap ((1-a) \cap \sigma(\bar a,\bar b)) = b$ which
  means $\pi(c) = b = c$.]

Now we can finish: as $\bbB_1 \subseteq \bbB'_2 \subseteq \bbB_2$ and
$\pi$ is a homomorphism from $\bbB'_2$ into $\bbB_1$ which is
a complete Boolean Algebra, we can extend
$\pi$ to $\pi^+$, a homomorphism from $\bbB_2$ into $\bbB_1$, see
\ref{x30}.  But $\pi$ is a projection hence so is $\pi^+$.  Clearly 
$(\pi^+)^{-1}\{1_{\bbB_1}\}$ includes $D$ and equality holds by
the assumption on the maximality of $D$ and we have proved the implication.
\medskip

\noindent
\underline{Clause (b) implies clause (a)}:

First, clearly $D$ is a filter of $\bbB_2$; also $a \in \bbB_1
\backslash \{1_{\bbB_1}\} \Rightarrow \pi(a) = a \ne 1_{\bbB_1}
\Rightarrow a \notin D$.

Toward contradiction assume $D_2$ is a filter on $\bbB_2,D \subsetneqq
D_2$ and $D_2 \cap \bbB_1 = \{1_{\bbB_1}\}$.  Choose $c_2 \in D_2
\backslash D$ and let $c_1 = \pi(c_2)$, consider the symmetric
difference, $c_1 \Delta c_2$ it is mapped by $\pi$ to $c_1 \Delta c_1 
= 0_{\bbB_2}$  hence $\pi(1_{\bbB_2} - (c_1 \Delta c_2)) = 1_{\bbB_2}
- \pi(c_1 \Delta c_2) = 1_{\bbB_2} - 0_{\bbB_2} = 1_{\bbB_2}$, so
$1_{\bbB_2} -(c_1 \Delta c_2) \in D$ so 
$c_1 = c_2 \mod D$, hence (recalling $D_1 \subseteq D_2$) we have
$c_1 = c_2 \mod D_2$ but $c_2 \in D_2$ hence $c_1 \in D_2$.  But
\mn
\begin{enumerate}
\item[$\bullet$]  $c_1 \in \bbB_1$ being $\pi(c_2)$
\sn
\item[$\bullet$]   $c_1 \ne 1_{\bbB_1}$ as $\pi(c_2) = c_1$ and
$c_2 \notin D$
\end{enumerate}
\mn
and recall
\mn
\begin{enumerate}
\item[$\bullet$]  $c_1 \in D_2$
\end{enumerate}
\mn
so $c_1$ contradicts $D_2 \cap \bbB_1 = \{1_{\bbB_1}\}$.  We comment
that for this direction we do not use the completeness of $\bbB_1$.

\noindent
1A) Now $\pi$ determines $D$ in the statement (b).  Also $D$
determines $\pi$ because if $\pi_1,\pi_2$ are projections from
$\bbB_2$ onto $\bbB_1$ such that $D = \{a \in \bbB_2:\pi_\ell(a) =
1_{\bbB_1}\}$ for $\ell=1,2$ and $\pi_1 \ne \pi_2$ let $a \in \bbB_2$
be such that $\pi_1(a) \ne \pi_2(a)$; then as in $(b) \Rightarrow (a)$
in the proof of part (1), $\pi_\ell(a) = a \mod D$ for $\ell=1,2$
hence $\pi_1(a) = \pi_2(a) \mod D$, but $\pi_1(a),\pi_2(a) \in 
\bbB_1$ and $D_2 \cap \bbB_1 = 1_{\bbB_1}$ and 
$D \cap \bbB_1 = 1_{\bbB_1}$ hence $\pi_1(a) = \pi_2(a)$, contradiction.

\noindent
2) Straightforward, but we elaborate; clearly $a \in D_1 \wedge b \in
D \Rightarrow \pi(a \cap b) = \pi(a) \cap \pi(b) = a \cap 0_{\bbB_1} =
a \ge 0_{\bbB_1}$ hence $a \in D_1 \wedge b \in D \Rightarrow a \cap b
> 0_{\bbB_2}$, hence least one such $D_2 \in \uf(\bbB_2)$.  For
uniqueness toward contradiction assume $\cE_1,\cE_2$ are from
$\uf(\bbB_2)$ and extend $D_1 \cup D$.  So necessarily there is $a \in
\cE_1 \backslash \cE_2$ so as above $a = \pi(a) \mod D$ but $D
\subseteq \cE_1 \cap \cE_2$ hence for $\ell=1,2$ we have $a = \pi(a)
\mod \cE_\ell$ so $a \in \cE_\ell \Leftrightarrow \pi(a) \in
\cE_\ell$.  But $\cE_1 \cap \bbB_1 = D_1 = \cE_2 \cap \bbB_1$ and
$\pi(a) \in \bbB_1$ has $\pi(a) \in \cE_1 \Leftrightarrow \pi(a) \in
\cE_2$.  By the last two sentences $a \in \cE_1 \Leftrightarrow a \in
\cE_2$ contradicting the choice of $a$.
\end{PROOF}

\begin{fact}  
\label{x33}
Assume $\bbB_1 \lessdot \bbB_2$ are complete Boolean Algebras,
$D_\ell \in \uf(\bbB_\ell)$ for $\ell=1,2$.  If $D$ is a maximal
filter on $\bbB_2$ disjoint to $\bbB_1 \backslash \{1_{\bbB_1}\}$ 
and $D \cup D_1
\subseteq D_2$ \then \, $\bbN^{\bbB_1}/D_1$ is an
initial segment of $\bbN^{\bbB_2}/D_2$.
\end{fact}

\begin{remark}  
\label{x36}
1) This is \cite[Ch.VI,Ex3.35]{Sh:a}.

\noindent
2) We can prove: if the homomorphism $\mathbf j:\bbB_2 
\rightarrow_{\text{onto}} \bbB_1$ maps $D_2 \in \uf(\bbB_2)$
   onto $D_1 \in \uf(\bbB_1)$ then $\bbN^{\bbB_1}/D_1$ is canonically
   isomorphic to an initial segment of $\bbN^{\bbB_2}/D_2$ as in \ref{x33}.
\end{remark}

\begin{PROOF}{\ref{x33}} 
The desired conclusion will follow by $(*)_3$ below:
\mn
\begin{enumerate}
\item[$(*)_1$]  If $\cI$ is a maximal antichain of $\bbB_1$ \then \,
$\{a/D:a \in \cI\}$ is a maximal antichain of $\bbB_2/D$.
\end{enumerate}
\mn
[Why?  First, 
\mn
\begin{itemize}
\item  $a \in \cI \Rightarrow a \in \bbB^+_1 \Rightarrow a/D
\in (\bbB_2/D)^+$
\sn
\item  if $a \ne b \in \cI$ then $\bbB_2 \models ``a \cap b =
  0_{\bbB_1}"$ hence $\bbB_2/D \models ``(a/D) \cap (b \cap D) =
  0_{\bbB_2/D}"$.
\end{itemize}
\mn
Hence, obviously $\cI^* := \{a/D:a \in \cI\}$ is an
antichain of $\bbB_2/D$.  Toward contradiction assume $\cI^*$ is not
maximal and let $c/D$ witness it.  By \ref{x23} there is $c' \in
\bbB_1$ such that $c=c' \mod D$ and so \wilog \, $c \in \bbB_1$.

As $c/D \ne 0/D$ necessarily $c \in \bbB^+_1$ and if $b \in \cI$ then
$(b/D) \cap (c/D) = 0/D$ hence $b \cap c = 0 \mod D$ but $b,c \in
\bbB_1$ hence $b \cap c=0$, so $c$ contradicts ``$\cI$ is a maximal
antichain of $\bbB_1$".]
\mn
\begin{enumerate}
\item[$(*)_2$]  If $f \in \bbN^{\bbB_2},c \in \bbB_1 \backslash D_1$
and $\TV(f > n) \cup c \in D$ for every $n$ \then \, $g \in
\bbN^{\bbB_1} \Rightarrow g/D_2 < f/D_2$.
\end{enumerate}
\mn
[Why?  If $g$ is a counter-example, then $\TV(f \le g)$ belongs to $D_2$ but
$1-c \in D_1 \subseteq D_2$ so $\TV(f \le g)-c$ belongs to $D_2$ 
hence to $D^+ := \{a \in \bbB_2:1-a \notin D\}$ since $D \subseteq
D_2$.   Let $\langle b_n:n < \omega\rangle$ represent 
$g$ as a member of $\bbN^{\bbB_1}$, 
then by $(*)_1, \, \langle
b_n/D:n < \omega\rangle$ is a maximal anti-chain of $\bbB_2/D$ hence
for some $n,\TV(f \le g) \cap b_n-c \in D^+$ but $\TV(f \le n) -c \ge
\TV(f \le g) \cap b_n-c$ hence $\TV(f \le n) -c \in D^+$,
contradiction to an assumption of $(*)_2$; so $(*)_2$ holds indeed.]
\mn
\begin{enumerate}
\item[$(*)_3$]  If $f \in \bbN^{\bbB_2},g \in \bbN^{\bbB_1}$ and
$f/D_2 \le g/D_2$ then for some $g' \in \bbN^{\bbB_1}/D$ we have
$f/D_2 = g'/D_2$.
\end{enumerate}
\mn
[Why?  Let $\langle a_n:n < \omega\rangle$ represent $f$ and let $a_{\ge
n} = \bigcup\limits_{k \ge n} a_k \in \bbB_2$.  If for 
some $b \in D_2$, we have $n < \omega \Rightarrow a_{\ge n} 
\cup(1-b) \in D$ then there is $f' \in \bbN^{\bbB_2}$ such that
$f'/D_2 = f/D_2$ and $n < \omega \Rightarrow \TV(f' \ge n) \in D$.
Now we apply $(*)_2$ with $f',0_{\bbB_1}$ here standing for $f,c$ there
and we get contradiction to ``$f/D_2 \le g/D_2$".  So we can assume
there is no such $b$.

Let $a'_n \in \bbB_1$ be such that $a_n = a'_n \mod D$ so possibly
 $a'_n = 0_{\bbB_1}$, such $a'_n$ exists by \ref{x23}(1A).  Clearly
 $\langle a'_n/D:n < \omega\rangle$ is an antichain of $\bbB_2/D$, so as
 $D \cap \bbB_1 = \{1_{\bbB_1}\}$ clearly $\langle a'_n:n < \omega\rangle$
is an antichain of $\bbB_1$.
\medskip  

\noindent
\underline{Case 1}:  $c := \bigcup\limits_{n} a'_n \notin D_1$.

First, we argue that for every $n \in \omega,\TV(f > n) \cup c \in D$ as 
otherwise there is some $n \in \omega$ such that $\TV(f \le n) -c \in
D^+$ hence for some $\ell \le n,a_\ell -c \in D^+$ hence (by the
choice of $a'_\ell$) we have $a'_\ell - c
\in D^+$, contradiction to the choice of $c$.  As $c \in \bbB_1$ we
see by the assumption of this case that $c \in \bbB_1 \backslash D_1$, hence 
by $(*)_2$ we get a contradiction to the assumption ``$f/D_2 \le
g/D_2$" of $(*)_3$. 
\medskip  

\noindent
\underline{Case 2}:  $c := \bigcup\limits_{n} a'_n \in D_1$ and $d =
\bigcup\limits_{n}(a_n \triangle a'_n) \notin D_2$.

As $D_2$ is an ultrafilter of $\bbB_2$ extending $D_1$, 
clearly $c' := c-d \in D_2$.  
We define $g' \in \bbN^{\bbB_1}$ as the function represented by
$\langle a'_n:n < \omega\rangle$ and $g'' \in \bbN^{\bbB_2}$ as the
function represented by 
$\langle a''_n:n < \omega \rangle$, where $a''_n$ is $a'_n \cap c'$ 
if $n > 0$ and $a'_n \cup(1-c')$ if $n=0$.  Easily 
$f/D_2 = g''/D_2$ because $f,g''$ ``agree" on $c'$ which belongs to
$D_2$ and is disjoint to $d$ and the choice of $d$; also 
$g''/D_2 = g'/D_2$ because $c' \in D_2$.  Together we are done.  
\medskip  

\noindent
\underline{Case 3}:  $c := \bigcup\limits_{n} a'_n \in D_1$ and $d =
\bigcup\limits_{n}(a_n \triangle a'_n) \in D_2$.

Let $d' \in \bbB_1$ be such that $d'/D = d/D$.  Let $d_1 :=
\bigcup\limits_{n}(a_n - a'_n)$ and $d_2 :=
\bigcup\limits_{n}(a'_n-a_n)$ hence $d=d_1 \cup d_2$.  Let $k <
\omega$, now modulo $D$ we have $d' \cap \bigcup\limits_{n \le k} a'_n = 
d \cap \bigcup\limits_{n \le k} a'_n =
\bigcup\limits^{2}_{\ell=1}(d_\ell \cap \bigcup\limits_{n \le k}
a'_n)$ and we shall deal separately with each term.  

First, $d_1 \cap \bigcup\limits_{n \le k} a'_n = 
\bigcup\limits_{\ell \le k} ((a_\ell - a'_\ell) \cap \bigcup\limits_{n \le k}
a'_n) \cup (\bigcup\limits_{\ell > k} (a_\ell - a'_\ell) 
\cap \bigcup\limits_{n \le k} a'_n)$.  Now the first term
$\bigcup\limits_{\ell \le k} ((a_\ell - a'_\ell) \cap \bigcup\limits_{n \le k}
a'_n)$ is equal $\mod D$ to $(\bigcup\limits_{n \le k} 0) \cap
\bigcup\limits_{n \le k} a'_k = 0_{\bbB_1}$, by the choice of the
$a'_\ell$.  Next, the second term in
the union, $(\bigcup\limits_{\ell > k} (a_\ell - a'_\ell) \cap
\bigcup\limits_{n \le k} a'_n)$ is modulo $D$ again by the choice of
the $a'_\ell$, equal to 
$(\bigcup\limits_{\ell > k} (a_\ell - a'_\ell)) \cap \bigcup\limits_{n \le k}
a_n$ which is zero as $\langle a_n:n < \omega\rangle$ is an antichain;
together by the previous sentences $d_1 \cap \bigcup\limits_{n \le k}
a'_n = 0_{\bbB_2} \mod D$.

Similarly $d_2 \cap \bigcup\limits_{n \le k} a'_n = 0_{\bbB_2} \mod D$
noting that $\langle a'_n:n < \omega\rangle$ is necessarily an
antichain of $\bbB_1$.  Hence $d' \cap \bigcup\limits_{n \le k} a'_n = d \cap
\bigcup\limits_{n \le k} a'_n = \bigcup\limits_{\ell =1}^{2} 
(d_\ell \cap \bigcup\limits_{n \le k} a'_n)) = 0_{\bbB_2} \cup
0_{\bbB_2} = 0_{\bbB_2} \mod D$.  But $d' \in \bbB_1$ and $a'_n \in
\bbB_1$ for every $n$ (and $D \cap \bbB_1 = \{1_{\bbB_1}\}$, of
course), hence $d' \cap \bigcup\limits_{n \le k} a'_n = 0_{\bbB_1}$.  However, 
as this holds for every $k$ and the choice of $c$ it follows that 
$d' \cap c =0$, but by the case first assumption of the present case
$c \in D_1 \subseteq D_2$ so $d' \notin D_2$, but by the 
case assumption $d'/D = d/D$ and $d \in D_2$ contradiction.
\end{PROOF}
\newpage

\section {Spectrum of the ultraproducts of finite cardinals}

\begin{definition}  
\label{k2}
Assume $D$ is an ultra-filter on $I$.

\noindent
1) Let $\upf(D)$ be the
spectrum of ultra-products $\mod D$ of finite cardinals, that is; 
$\{\prod\limits_{i \in I} n_i/D:n_i \in \bbN$ for 
$i \in I$ and $\prod\limits_{i \in I} n_i/D$ is infinite$\}$.

\noindent
2) For $\lambda \in \upf(D)$ let $A_{D,\lambda} = \{a:a \in \bbN^I/D$
and the set $\{b \in \bbN^I/D:\bbN^I/D \models ``b < a"\}$ has cardinality
$\lambda\}$; we consider it as a linearly ordered set by the order
inherited from $\bbN^I/D$.

\noindent
3) Let $\spec_1(D) = \{(\lambda,\theta,\partial):\lambda \in \upf(D)$
   and $A_{D,\lambda}$ has cofinality $\partial$ and co-initiality
   $\theta\}$.

\noindent
4) Let $\spec_2(D)$ be the set of triples 
$(\lambda,\theta,\partial)$ such that:
\mn
\begin{enumerate}
\item[$(a)$]  $\lambda \in \upf(D)$
\sn
\item[$(b)$]  $(\alpha) \quad$ if $\lambda < \max(\upf(D))$ then
$A_{D,\lambda}$ has cofinality $\partial$
\sn
\item[${{}}$]  $(\beta) \quad$ if $\lambda = \max(\upf(D))$ then
$\partial = 0$ (or $*$)
\sn
\item[$(c)$]  $\theta$ is the co-initiality of $A_{D,\lambda}$.
\end{enumerate}
\mn
5) For $D$ an ultrafilter on a complete Boolean Algebra $\bbB$ we define the
above similarly considering $\bbN^{\bbB}/D$ instead $\bbN^J/D$ but in
clause (b), $\partial$ is the cofinality of $A_{D,\lambda}$ in all cases.
\end{definition}

\begin{definition}
\label{f2}
Let $K_\alpha$ be the class of objects $\mathbf k$ consisting of:
\mn
\begin{enumerate}
\item[$(a)$]  $\bbB_\beta$ is a Boolean Algebra for $\beta \le \alpha$
\sn
\item[$(b)$]  $\langle \bbB_\beta:\beta \le \alpha \rangle$ is increasing
\sn
\item[$(c)$]  $\bbB_\beta$ is complete for $\beta < \alpha,\bbB_0$ is trivial
\sn
\item[$(d)$]  $\bbB_\beta \lessdot \bbB_\gamma$ if $\beta < \gamma \le
\alpha$ and $\cup\{\bbB_{\beta_1}:\beta_1 < \gamma\} \lessdot
\bbB_\gamma$ for limit $\gamma \le \alpha$
\sn
\item[$(e)$]  $D_\beta$ is a filter on $\bbB_\alpha$ such that
$\bbB_\beta \cap D_\beta = \{1_{\bbB_\beta}\}$
\sn
\item[$(f)$]  $D_\beta$ is maximal under clause $(e)$, so $D_0$ is an
ultrafilter and $D_\alpha = \{1_{\bbB_\alpha}\}$
\sn
\item[$(g)$]  $\langle D_\beta:\beta \le \alpha\rangle$ is 
$\subseteq$-decreasing.
\end{enumerate}
\end{definition}

\begin{definition}
\label{f4}
1) Above let $\bbB[\mathbf k] = \bbB_{\mathbf k} = \bbB_\alpha,\bbB[\mathbf
   k,\beta] = \bbB_{\mathbf k,\beta} = \bbB_\beta,\bar{\bbB}_{\mathbf k} = \langle
   \bbB_{\mathbf k,\beta}:\beta \le \alpha\rangle,D_{\mathbf k,\beta} =
   D_\beta,D_{\mathbf k} = D_{\mathbf k,0},\ell g(\mathbf k) = \alpha_{\mathbf k} =
   \alpha(\mathbf k) = \alpha$.

\noindent
1A) Let $K^{\com}_\alpha$ be the class of $\mathbf k \in K_\alpha$ such
that $\bbB_{\mathbf k}$ is a complete Boolean Algebra.

\noindent
2) Assume $\kappa > \aleph_0$ is regular.
Let $K^{\cc(\kappa),1}_\alpha$ be the class of $\mathbf k \in
   K_\alpha$ such that $\bbB_\alpha$ satisfies the $\kappa$-c.c.

\noindent
3) Let $K^{\cc(\kappa),2}_\alpha$ be the class of $\mathbf k \in
   K^{\cc(\kappa),1}_\alpha$ such that:
\mn
\begin{enumerate}
\item[$\bullet$]   $\bbB_{\mathbf k}$ is complete; recall that for
  every $\beta < \alpha,\bbB_{\beta}$ is complete 
\sn
\item[$\bullet$]  if $\delta \le \alpha$ has cofinality $\ge \kappa$
\then \, $\bbB_{\mathbf k,\delta} = \bigcup\limits_{\beta < \delta} 
\bbB_{\mathbf k,\beta}$
\sn 
\item[$\bullet$]  if $\delta \le \alpha$ is limit of cofinality $<
\kappa$, \then \,
$\bbB_{\mathbf k,\delta}$ is the completion of $\bigcup\limits_{\beta <
\delta} \bbB_{\mathbf k,\beta}$.
\end{enumerate}
\mn
3A) We may omit $\kappa$ when $\kappa = \aleph_1$ so $K^{\cc,\iota}_\alpha =
K^{\cc(\aleph_1),\iota}_\alpha$; if we omit $\iota$ we mean 1.

\noindent
4) Let $K = \bigcup\limits_{\alpha} K_\alpha$ and
$K^{\cc(\kappa),\iota} = \cup\{K^{\cc(\kappa),\iota}_\alpha:\alpha$ is
an ordinal$\}$ so $K^{\cc} = \bigcup\limits_{\alpha} K^{\cc}_\alpha$.

\noindent
5) We say $\bbB$ is above $\bar{\bbB}_{\mathbf k}$ when $\bbB_{\mathbf k} 
\subseteq \bbB$ and $\bbB_{\mathbf k,\beta} \lessdot \bbB$
 for $\beta < \alpha_{\mathbf k}$.

\noindent
6) $K^{\fr(\kappa)}_\alpha$ is the class of $\mathbf f$ consisting of:
\mn
\begin{enumerate} 
\item[$(a)$]  $\mathbf k_{\mathbf f} = (\bar{\bbB},\bar D)$ as in \ref{f2}
\sn
\item[$(b)$]  $\bar \xi = \langle \xi_\gamma:\gamma \le \alpha\rangle$ and
$\bar x = \langle x_{\beta,\zeta,i}:i < \kappa,\beta < \alpha,\zeta <
\xi_\beta\rangle,x_{\beta,\zeta,i} \in \bbB_{\mathbf k}$ are such 
that $\bar x$ is free except
that $\beta < \xi_\alpha \wedge i < j < \kappa \Rightarrow 
x_{\beta,\zeta,i} \cap x_{\beta,\zeta,j} = 0$
\sn
\item[$(c)$]  the sub-algebra which $\langle x_{\beta,\zeta,i}:\zeta <
  \xi_\gamma,i < \kappa\rangle$ generates is dense in $\bbB_{\mathbf
  k,\gamma}$
\sn
\item[$(d)$]  so $\bar\xi_{\mathbf f} = \bar\xi,\bar x_{\mathbf f} = \bar
  x,\bar{\bbB}_{\mathbf f} = \bbB_{\mathbf k}$, etc.
\end{enumerate}
\mn
7) Let ${}_* K_\alpha$ be defined like $K_\alpha$ in \ref{f2} omitting
clause (d) of \ref{f2}, and define ${}_* K$, as above; not really 
needed here but we may comment.
\end{definition}

\begin{definition}
\label{f8}
1) If $\beta \le \gamma$ and $\mathbf m \in K_\gamma$ then $\mathbf k = \mathbf m
   \rest \beta$ is the unique $\mathbf k \in K_\beta$ such that
   $\bbB_{\mathbf k} = \bbB_{\bfm,\beta},
\bbB_{\mathbf k,\alpha} = \bbB_{\mathbf m,\alpha},D_{\mathbf k,\alpha} =
   D_{\mathbf m,\alpha} \cap \bbB_{\mathbf k}$ for $\alpha \le \beta$. 

\noindent
1A) If $\mathbf k \in K_\alpha$ and $\beta < \alpha$ \then \,
$\pi_{\mathbf k,\beta}$ is the unique projection from $\bbB_{\mathbf k}$
onto $\bbB_{\mathbf k,\beta}$ such that 
$\pi^{-1}_{\mathbf k,\beta}\{1_{\bbB_{\mathbf k,\beta}}\} = 
D_{\mathbf k,\beta}$ recalling \ref{x23}; let
$\pi_{\mathbf k,\alpha} = \id_{\bbB_{\mathbf k,\alpha}}$ and if $\gamma
\le \beta \le \alpha$ then $\pi_{\mathbf k,\beta,\gamma} = \pi_{\mathbf
  k,\gamma} \rest \bbB_{\mathbf k,\beta}$.

\noindent
2) We define the following two-place relations on $K$:
\mn
\begin{enumerate}
\item[$(A)$]  $\mathbf k \le^{\at}_K \mathbf m$ where $\at$ stands 
for atomic \Iff \,:
\sn
\begin{enumerate}
\item[$(a)$]  $\alpha_{\mathbf k} = \alpha_{\mathbf m}$
\sn
\item[$(b)$]  $\bbB_{\mathbf k,\beta} = \bbB_{\mathbf m,\beta}$ for $\beta
< \alpha_{\mathbf k}$
\sn
\item[$(c)$]  $\bbB_{\mathbf k,\alpha(\mathbf k)} \lessdot
\bbB_{\mathbf m,\alpha(\mathbf m)}$
\sn
\item[$(d)$]  $D_{\mathbf k,\beta} \subseteq D_{\mathbf m,\beta}$ for $\beta
\le \alpha_{\mathbf k}$.
\end{enumerate}
\sn
\item[$(B)$]  $\mathbf k \le^{\ver}_K \mathbf m$, where $\ver$ stands 
for vertical \Iff \,
\sn
\begin{enumerate}
\item[$(a)$]  $\alpha_{\mathbf k} \le \alpha_{\mathbf m}$
\sn
\item[$(b)$]  $\mathbf k \le^{\at}_K (\mathbf m \rest \alpha_{\mathbf k})$
\end{enumerate}
\sn
\item[$(C)$]  $\mathbf k \le^{\hor}_K \mathbf m$, where $\hor$ stands 
for horizontal \Iff \,
\sn
\begin{enumerate}
\item[$(a)$]  $\alpha_{\mathbf k} = \alpha_{\mathbf m} = \alpha$
\sn
\item[$(b)$]  $\bbB_{\mathbf k,\beta} \lessdot \bbB_{\mathbf m,\beta}$ for $\beta
\le \alpha$
\sn
\item[$(c)$]  $D_{\mathbf k,\beta} \subseteq D_{\mathbf m,\beta}$ for
$\beta \le \alpha$
\end{enumerate}
\sn
\item[$(D)$]
\begin{enumerate}
\item[$(a)$]  $\mathbf f_1 \le^{\fr(\kappa)}_K \mathbf f_2$ iff
$\mathbf f_\ell \in K^{\fr(\kappa)}_{\alpha_\ell}$ and $\mathbf k_{\mathbf
  f_1} \le^{\ver}_K \mathbf k_{\mathbf f_2}$ and $\bar x_{\mathbf f_1}
\trianglelefteq \bar x_{\mathbf f_2}$ which means:
\sn
\begin{itemize}
\item  $\beta < \alpha_1 \Rightarrow \bar x_{\mathbf f_1,\beta} = \bar
  x_{\mathbf f_2,\beta}$ and $\beta = \alpha_1 \Rightarrow \xi_{\mathbf
    f_1,\beta} \le \xi_{\mathbf f_2,\beta} \wedge \bar x_{\mathbf
    f_1,\beta} = \bar x_{\mathbf f_2,\beta} \rest \xi_{\mathbf f_2,\beta}$
\end{itemize}
\sn
\item[$(b)$]  $\le^{\at-\fr(\kappa)}_K$ is defined similarly
\end{enumerate}
\sn
\item[$(E)$]  $\mathbf k \le_K^{\wa} \mathbf m$ where $\wa$ stands 
for weakly atomic \Iff \,
\sn
\item[${{}}$]  $(a),(b),(d) \quad$ as in Clause (A)
\sn
\item[${{}}$]  $(c) \quad \bbB_{\mathbf k,\alpha(\mathbf k)} \subseteq
  \bbB_{\mathbf m,\alpha(\mathbf m)}$.
\end{enumerate}
\end{definition}

\begin{remark}
\label{f3}
Note that for the present work it is not a loss to use exclusively
c.c.c. Boolean Algebras; moreover ones which have a dense subalgebra which
is free.  So using only free
Boolean Algebras or their completion, i.e. 
$(K^{\fr(\aleph_1)}_\alpha,\le^{\fr(\kappa)}_K)$; so we are 
giving for $\bbB$ a set of generators (and the orders respect this).
\end{remark}

\begin{observation}
\label{f10}
The relations $\le^{\at}_K,\le^{\wa}_K$ and 
$\le^{\ver}_K$ and $\le^{\hor}_K$ (the
last one is not used) are partial orders on $K$.
\end{observation}

\noindent
We need various claims on extending members of $K$, existence of upper
bounds to an increasing sequence and amalgamation.
\begin{claim}
\label{f13}
Let $\delta$ be a limit ordinal.

\noindent
1) If $\langle \mathbf k_i:i < \delta\rangle$ is
$\le^{\at}_K$-increasing \then \, it has a $\le^{\at}_K$-lub 
$\mathbf k_\delta$, the union naturally defined so $|\bbB_{\mathbf k_\delta}|
\le \Sigma\{|\bbB_{\mathbf k_i}|:i < \delta\}$.

\noindent
1A) Like part (1) for $\le^{\wa}_K$.

\noindent
2) If $\langle \mathbf k_i:i < \delta\rangle$ is a $\le^{\ver}_K$-increasing
   sequence, \then \, it has a $\le^{\ver}_K$-upper bound $\mathbf k =
   \mathbf k_\delta$ which is the union which means:
\mn
\begin{enumerate}
\item[(a)]  $\ell g(\mathbf k) = \cup\{\ell g(\mathbf k_i):i < \delta\}$
call it $\alpha$
\sn
\item[(b)]  if $\beta < \alpha$ then $\bbB_{\mathbf k,\beta} =
\bbB_{\mathbf k_i,\beta}$ for every large enough $i$
\sn
\item[(c)]
\begin{enumerate}
\item[$(\alpha)$]  if $\langle \ell g(\mathbf k_i):i <
\delta\rangle$ is eventually constant (so $\ell g(\mathbf k_i) = \alpha$
for every $i <\delta$ large enough) \then \,
\sn
\begin{itemize}
\item  $\bbB_{\mathbf k,\alpha} =
\cup\{\bbB_{\mathbf k_i,\alpha}:i < \delta$ is such that $\ell g(\mathbf k_i)
= \alpha\}$
\sn
\item  $D_{\mathbf k,\alpha} = \{1_{\bbB_{\mathbf k,\alpha}}\}$, redundant
\end{itemize}
\sn
\item[$(\beta)$]  if $\langle \ell g(\mathbf k_i):i <
\delta\rangle$ is not eventually constant \then
\sn
\begin{itemize}
\item  $\bbB_{\mathbf k,\alpha} =
\cup\{\bbB_{\mathbf k_i,\ell g(\mathbf k_i)}:i < \delta\}$
\sn
\item   $D_{\mathbf k,\alpha} = \{1_{\bbB_{\mathbf k,\alpha}}\}$, redundant
\end{itemize}
\end{enumerate}
\sn
\item[(d)]  if $\beta < \alpha$ then $D_{\mathbf k,\beta} = 
\cup\{D_{\mathbf k_i,\beta}:i < \delta$ is such that 
$\beta \le \ell g(\mathbf k_i)\}$.
\end{enumerate}
\mn
3) In part (2) if $\mathbf k_i \in K^{\cc(\kappa),2},\cf(\delta) \ge
\kappa$ and the sequence $\langle \ell g(\mathbf k_i):i < \delta\rangle$ 
is not eventually constant \then \, $\bbB_{\mathbf k}$ is complete and
$\upf(D_{\mathbf k}) = \cup\{\upf(D_{\mathbf k_i}):i < \delta\}$.

\noindent
4) Similarly for the ${}_* K$ version.
\end{claim}

\begin{PROOF}{\ref{f13}}
Straightforward (concerning part (3) note that recalling $\cf(\delta)
\ge \kappa$ we have $\bbN^{\bbB(\mathbf k)} = \cup\{\bbN^{\bbB(\mathbf
  k_i)}:i < \delta\}$.
\end{PROOF}

\begin{claim}
\label{f16}
1) If $\mathbf k \in K_\alpha$ and $\bbB$ is above $\bar{\bbB}_{\mathbf k}$
(i.e. $\bbB_{\mathbf k,\alpha} \subseteq \bbB$ and $\beta < \alpha \Rightarrow
\bbB_{\mathbf k,\beta} \lessdot \bbB$) \then \, there is $\mathbf m \in
K_\alpha$ such that $\mathbf k \le^{\wa}_K \mathbf m$ and $\bbB_{\mathbf
m,\alpha} = \bbB$.

\noindent
2) If $\mathbf k \in K_\alpha$ and $\bbB \subseteq \bbB_{\mathbf k}$ and
   $\beta < \alpha \Rightarrow \bbB_{\mathbf k,\beta} \subseteq \bbB$
   (hence $\bbB_{\mathbf k,\beta} \lessdot \bbB$) \then \, there is
   $\mathbf m \in K_\alpha$ such that $\mathbf m \le^{\wa}_K \mathbf k$ and
   $\bbB_{\mathbf m} = \bbB$. 

\noindent
3) If $\mathbf k \in K_\alpha,\bbB_{\mathbf k} \lessdot \bbB$ \then \, for
   some $\mathbf m \in K_\alpha,\mathbf k \le^{\at}_K \mathbf m,\bbB_{\mathbf
   m} = \bbB$

\noindent
4) In part (2) if $\bbB \lessdot \bbB_{\mathbf k}$ \then \, we can add
   $\mathbf m \le^{\at}_K \mathbf k$.

\noindent
5) If $\mathbf k \le^{\wa}_K \mathbf m$ and $\bbB_{\mathbf k} \lessdot
   \bbB_{\mathbf m}$ \then \, $\mathbf k \le^{\at}_K \mathbf m$.
\end{claim}

\begin{PROOF}{\ref{f16}}
1) For simplicity we concentrate on the main case: all the
$\bbB_\beta$ and $\bbB$ satisfies the c.c.c.
By \ref{f13}(1A) \wilog \, $\bbB$ is generated by $\bbB_{\mathbf k,\alpha}
\cup \{a\}$ where $a \notin \bbB_{\mathbf k,\alpha}$.  So $\bbB_{\mathbf m}$
is uniquely defined and as required in Definition \ref{f2}, but we
have to define the $D_{\mathbf m,\beta}$'s and, of course, let
$D_{\mathbf m,\alpha} = \{1_{\bbB_{\mathbf m}}\}$.
\medskip

\noindent
\underline{Case 1}:  $\alpha = 0$

This is trivial.
\medskip

\noindent
\underline{Case 2}:  $\alpha = \beta +1$

As $\bbB_{\mathbf k,\beta} \lessdot \bbB$ and $\pi_{\mathbf k,\beta}$ is a
projection from $\bbB_{\mathbf k,\alpha}$ onto $\bbB_{\mathbf k,\beta}$
and $[B_{\mathbf k,\alpha} \subseteq \bbB$ and $\bbB_{\mathbf k,\beta}$ is
complete by \ref{x30}(1) there is a projection $\pi$
from $\bbB$ onto $\bbB_{\mathbf k,\beta}$ extending 
$\pi_{\mathbf k,\beta}$.  Now for $\gamma < \alpha$ let $D_{\mathbf m,\gamma}$
be the filter on $\bbB$ generated by $D_{\mathbf k,\gamma} \cup 
\{(a \triangle \pi(a)\}$.
\medskip

\noindent
\underline{Case 3}:  $\alpha$ is a limit ordinal, $\bbB$ is
c.c.c. (the main case) and $\alpha$ is of cofinality $> \aleph_0$

In this case, $\bbB' = \bigcup\limits_{\gamma < \alpha} \bbB_{\mathbf
  k,\gamma}$ is complete and $\lessdot \bbB$ so we can continue as in
Case 2 using $\bbB'$ instead $\bbB_{\mathbf k,\beta}$.
\medskip

\noindent
\underline{Case 4}:  $\cf(\alpha) = \aleph_0$

Let $\alpha = \bigcup\limits_{n} \alpha_n,\alpha_n < \alpha_{n+1}$.
For $\beta < \alpha$ let $\pi_\beta$ be the projection of $\bbB_{\mathbf
  k}$ onto $\bbB_{\mathbf k,\beta}$ which maps $D_{\mathbf k,\beta}$ onto
$1_{\bbB_{\mathbf k,\beta}}$.  Let $\Pi_\beta$ be the set of
homomorphisms from $\bbB$ into $\bbB_{\mathbf k,\beta}$ extending
$\pi_\beta = \pi_{\bfk,\beta}$, so not empty hence 
(recalling $\bbB_{\mathbf k,\beta}$ is complete)
there are $b_\beta \le c_\beta$ from 
$\bbB_{\mathbf k,\beta}$ such that $\{\pi(a):\pi \in\Pi_\beta\}$ is
$\{a' \in \bbB_{\mathbf k,\beta}:b_\beta \le a' \le c_\beta\}$.  Clearly
$\gamma < \beta < \alpha \Rightarrow b_\gamma \le b_\beta \le c_\beta
\le c_\gamma$ in $\bbB_{\mathbf k,\beta}$.

Now by induction on $\zeta < (\|\bbB\|^{|\alpha|})^+$ we define
$\langle (b_{\beta,\zeta},c_{\beta,\zeta}):\beta < \alpha\rangle$ such
that:
\mn
\begin{enumerate}
\item[$(*)_\zeta$]
\begin{enumerate}
\item[(a)]  $\bbB_{\mathbf k,\beta} \models ``b_{\beta,\zeta} \le
  c_{\beta,\zeta}"$ for $\beta < \alpha$
\sn
\item[(b)]  if $\gamma < \beta < \alpha$ then 
$\bbB_{\mathbf k,\beta} \models ``b_{\gamma,\zeta} \le b_{\beta,\zeta}
\le c_{\beta,\zeta} \le c_{\gamma,\zeta}"$
\sn
\item[(c)]  if $\varp < \zeta$ and $\beta < \alpha$ then 
$\bbB_{\mathbf k,\beta} \models ``b_{\beta,\varp} \le b_{\beta,\zeta}
  \le c_{\beta,\zeta} \le c_{\beta,\varp}"$.
\end{enumerate}
\end{enumerate}
\medskip

\noindent
\underline{Subcase 4A}:  For $\zeta = 0$ let
$(b_{\beta,\zeta},c_{\beta,\zeta}) = (b_\beta,c_\beta)$, so clause
(a),(b) holds as said above and clause (c) is empty.
\medskip

\noindent
\underline{Subcase 4B}:  $\zeta$ is a limit ordinal

Let for $\beta < \alpha$
\mn
\begin{itemize}
\item  $b_{\beta,\zeta} = \cup\{b_{\gamma,\varp}:
\varp < \zeta\}$ in $\bbB_{\mathbf k,\beta}$
\sn
\item  $c_{\beta,\zeta} = \cap\{b_{\gamma,\varp}:\varp < \zeta\}$.
\end{itemize}
\mn
They are well defined because $\bbB_{\mathbf k,\beta}$ is a complete
Boolean Algebra and it is easy to check that (a),(b),(c) hold.
\medskip

\noindent
\underline{Subcase 4C}:  $\zeta = \varp + 1$

Let $b_{\beta,\zeta} = \cup\{\pi_{\mathbf k,\gamma,\beta}
(b_{\gamma,\varp}):\gamma \in (\beta,\alpha)\},c_{\beta,\zeta} =
\cap\{\pi_{\mathbf k,\gamma,\beta}(c_{\gamma,\varp}):\gamma \in 
(\beta,\alpha)\}$.

Now check.

Having carried the induction, by $(*)_\zeta$ for $\zeta <
(\|\bbB\|^{|\alpha|})^+$ for some $\zeta_* <
(\|\bbB\|^{|\alpha|})^+,\langle
(b_{\beta,\zeta},c_{\beta,\zeta}):\beta < \alpha\rangle$ is the same
for all $\zeta \ge \zeta_*$ and let $a_\beta =
b_{\beta,\zeta_*}$ for $\beta < \alpha$.

Easily $\gamma < \beta < \alpha \Rightarrow \pi_\beta(a_\beta) =
a_\gamma$ and let $D_{\mathbf m,\beta}$ be the filter of $\bbB$ generated by
$D_{\mathbf k,\beta} \cup \{(a \triangle a_\beta)\}$.

\noindent
2)  Easy.

\noindent
3) By (1).

\noindent
4),5)  Should be easy.
\end{PROOF}

\begin{claim}
\label{f17}
1) If $\mathbf k \in K_\alpha,\langle \beta_i:i \le i(*)\rangle$ is
increasing with $\beta_{i(*)} = \alpha$ \then \, there is one and only
one $\mathbf m \in K_{i(*)}$ such that $(\bbB_{\mathbf m,i},D_{\mathbf m,i})
= (\bbB_{\mathbf k,\beta_i},D_{\mathbf m,\beta_i})$ for $i \le i(*)$.

\noindent
2) Above if $\mathbf m \le^{\at}_K \mathbf m_1$ \then \, there is $\mathbf
k_1$ such that $\mathbf k \le^{\at}_K \mathbf k_1$ and $\bbB_{\mathbf k_1} =
\bbB_{\mathbf m_1}$ and $D_{\mathbf k_1,\beta_i} = D_{\mathbf m_1,i}$ for $i
\le i(*)$.

\noindent
2A) Similarly for $\le^{\wa}_K$.

\noindent
3) Above if $\mathbf m \le^{\ver}_K \mathbf m_1 \in K_{i(*) + j(*)}$ \then \,
there is $\mathbf k_1 \in K_{\alpha +j(*)}$ such that $\mathbf k
\le^{\ver}_K \mathbf k_1,\bbB_{\mathbf k_1} = \bbB_{\mathbf m_1},\bbB_{\mathbf
  k_1,\alpha +j} = \bbB_{\mathbf m_1,i(*) + j},D_{\mathbf k_1,\alpha +j}
= D_{\mathbf m_1,i(*) + j},D_{\mathbf k_1,\beta_i} = D_{\mathbf m_1,i}$
for $j < j(*),i \le i(*)$.

\noindent
4) If $\mathbf k \in K_\alpha$ and $\beta_0 \le \beta_1 \le \beta_2 \le
   \alpha$ \then \, $\pi_{\mathbf k,\beta_2,\beta_0} = \pi_{\mathbf
   k,\beta_1,\beta_0} \circ \pi_{\mathbf k,\beta_2,\beta_1}$.
\end{claim}

\begin{PROOF}{\ref{f17}}
Straightforward.
\end{PROOF}
 
\begin{conclusion}
\label{f19}
1) If $\mathbf k \in K_\alpha$ \then \, there is $\mathbf m$ such that
   $\mathbf k \le^{\at}_K \mathbf m$ and $\bbB_{\mathbf m,\alpha}$ is the
   completion of $\bbB_{\mathbf k,\alpha}$ so $\mathbf m \in K^{\com}_\alpha$ and
   $\mathbf k \le^{\ver}_K \mathbf m$; so if $\mathbf k \in
K^{\cc(\kappa)}_\alpha$ then $\mathbf m \in K^{\cc(\kappa)}_\alpha \cap
K^{\com}$.

\noindent
2) If $\alpha < \beta,\mathbf k \in K_\alpha,\mathbf n \in K_\beta$ and
$\mathbf k \le^{\ver}_K \mathbf n$, \then \, for some $\mathbf m$ we have
$\mathbf k \le^{\ver}_K \mathbf m \le^{\ver}_K \mathbf n$ and $\bbB_{\mathbf
  m}$ is the completion of $\bbB_{\mathbf k}$ inside $\bbB_{\mathbf n}$;
so if $\mathbf n \in K^{\cc(\kappa)}$ then $\mathbf m \in K^{\cc(\kappa)}$.
\end{conclusion}

\begin{PROOF}{\ref{f19}}
1) By \ref{f16}.

\noindent
2) Check the definitions.
\end{PROOF}

\begin{claim}
\label{f22}
There is $\mathbf m$ such that $\mathbf k \le^{\wa}_K \mathbf m,
\bbB_{\mathbf m} = \bbB$ and $Y \subseteq D_{\mathbf m}$ and $\bbB_{\mathbf
  k} \lessdot \bbB \Rightarrow \mathbf k \le^{\at}_K \mathbf m$ \when \,:
\mn
\begin{enumerate}
\item[(a)]  $\mathbf k \in K_\alpha$
\sn
\item[(b)]  $\bbB$ is a Boolean Algebra, $Y \subseteq \bbB$
\sn
\item[(c)]
\begin{enumerate}
\item[$(\alpha)$]  $\bbB_{\mathbf k} \subseteq \bbB$
\sn
\item[$(\beta)$]  $\bbB_{\mathbf k,\beta} \lessdot \bbB$ for $\beta <
\alpha_{\mathbf k}$
\end{enumerate}
\sn
\item[(d)]   if $\beta < \alpha_{\mathbf k}$ \then \, for some
$X_\beta$ we have:
\sn
\begin{enumerate}
\item[$(\alpha)$]  $X_\beta \subseteq Y$
\sn
\item[$(\beta)$]  $X_\beta$ is a downward directed subset of $\bbB$
\sn
\item[$(\gamma)$]   if $x \in X_\beta$ and $b \in
D_{\mathbf k,\beta}$ \then \, $x \cap b$ is not disjoint to any
$a \in \bbB^+_{\mathbf k,\beta}$
\sn
\item[$(\delta)$]  if $y \in Y$ \then \, for some $b \in D_{\mathbf
  k,\beta}$ and $x \in X_\beta$ we\footnote{Does this contradict
  $(d)(\gamma)$?  No, as $\mathbf D_{\mathbf k,\beta}$ is disjoint to
$\bbB_{\mathbf k,\beta} \backslash \{1_{\bbB_{\mathbf k,\beta}}\}$.}
 have $b \cap x \le y$.
\end{enumerate}
\end{enumerate}
\end{claim}

\begin{PROOF}{\ref{f22}}
By \ref{f16} \wilog \, $\bbB$ is generated by $\bbB_{\mathbf k} \cup Y$
and let $\alpha = \alpha_{\mathbf k}$ and let 
$X_\beta$ be as in clause (d) in the claim for
 $\beta < \alpha$ and define $\mathbf m \in K_\alpha$ as follows:
\mn
\begin{enumerate}
\item[$\bullet$]  $D_{\mathbf m,\alpha} = \{1_{\bbB}\},
\bbB_{\mathbf m,\alpha} = \bbB$ 
\sn
\item[$\bullet$]  for $\beta < \alpha$ let
$\bbB_{\mathbf m,\beta} = \bbB_{\mathbf k,\beta}$ and 
$D_{\mathbf m,\beta}$ be the filter on $\bbB_{\mathbf m,\alpha}$
  generated by $D_{\mathbf k,\beta} \cup X_\beta$.
\end{enumerate}
\mn
The point is to check $\mathbf m \in K_\alpha$ as then $\mathbf k
\le^{\at}_K \mathbf m$ and $Y \subseteq D_{\mathbf m}$ are obvious, also
$\bar{\bbB}_{\mathbf m}$ is as required and $D_{\mathbf m,\beta}$ a filter
on $\bbB_{\mathbf m,\alpha}$ including $D_{\mathbf k,\beta}$ are obvious.

So proving $(*)_1,(*)_2,(*)_3$ below will suffice
\mn
\begin{enumerate}
\item[$(*)_1$]  if $\beta < \gamma < \alpha$ then $D_{\mathbf m,\gamma}
  \subseteq D_{\mathbf m,\beta}$.
\end{enumerate}
\mn
[Why?  If $a \in D_{\mathbf m,\gamma}$ then by the choice of 
$D_{\mathbf m,\gamma}$ (recalling $D_{\mathbf k,\gamma}$ is downward
  directed being a filter and $X_\gamma$ is downward directed by its
  choice (i.e. Clause $(d)(\beta)$ of the claim) for 
some $b \in D_{\mathbf k,\gamma}$ and 
$x \in X_\gamma$ we have $b \cap x \le a$.  So by $(d)(\alpha)$
  applied to $\gamma$ we have $x \in Y$ hence 
by $(d)(\delta)$ applied to $\beta$ for
  some $b_1 \in D_{\mathbf k,\beta}$ and $x_1 \in X_\beta$ we have $b_1
 \cap x_1 \le x$ hence $(b \cap b_1) \cap x_1 = b \cap (b_1
  \cap x_1) \le b \cap x \le a$ but $b \in D_{\mathbf k,\gamma}
  \subseteq D_{\mathbf k,\beta},b_1 \in D_{\mathbf k,\beta}$ hence $b \cap
  b_1 \in D_{\mathbf k,\beta}$ and $x_1 \in X_\beta$ hence $a \in
  D_{\mathbf m,\beta}$ by the choice of $D_{\mathbf m,\beta}$.]
\mn
\begin{enumerate}
\item[$(*)_2$]   $D_{\mathbf m,\beta}$ is a filter on 
$\bbB_{\mathbf m,\alpha} = \bbB$ disjoint to $\bbB_{\mathbf m,\beta}
\backslash \{1_{\bbB_{\mathbf k}}\} = \bbB_{\mathbf k,\beta} \backslash
\{1_{\bbB_{\mathbf k}}\}$.
\end{enumerate}
\mn
[Why?  By the definition of $D_{\mathbf m,\beta}$ and clause
  $(d)(\gamma)$.]
\mn
\begin{enumerate}
\item[$(*)_3$]   if $\beta < \alpha$ then
$D_{\mathbf m,\beta}$ is a maximal filter of 
$\bbB_{\mathbf m}$ disjoint to $\bbB_{\mathbf k,\beta} \backslash
  \{1_{\mathbf m,\beta}\} = \bbB_{\mathbf m,\beta} \backslash
  \{1_{\bbB_{\mathbf m,\beta}}\}$.
\end{enumerate}
\mn
Why?  If $b \in \bbB = \bbB_{\mathbf m,\alpha}$ then for some Boolean terms
$\sigma(y_0,y_1,\dotsc,z_0,z_1,\ldots)$ and $a_0,a_1,\dotsc,\in 
\bbB_{\mathbf k}$ and $x_0,x_1,\ldots \in
Y$ we have $b = \sigma(a_0,a_1,\dotsc,x_0,x_1,\ldots)$ hence modulo
 the filter $D_{\mathbf m,\beta},b$ is equal to
$\sigma(a_0,a_1,\dotsc,1_{\bbB_{\mathbf k,\beta}},1_{\bbB_{\mathbf k,\beta}},
\ldots)$.  But for each $a_\ell$ there is $c_\ell \in \bbB_{\mathbf
  k,\beta}$ such that $a_\ell = c_\ell \mod D_{\mathbf k,\beta}$ hence
$b$ is equal to $\sigma(c_0,c_1,\dotsc,1_{\bbB_{\mathbf
    k,\beta}},1_{\bbB_{\mathbf k,\beta}},\ldots)$ which belongs 
to $\bbB_{\mathbf k,\beta}$.

As this holds for any $b \in \bbB$ we are easily done.
\end{PROOF}

\begin{definition}
\label{f24}
1) We say $\mathbf k$ is reasonable in $\alpha$ \when \, $\alpha +1 \le
\alpha_{\mathbf k}$ (so $\bbB_{\mathbf k,\alpha}$ is complete) and there is a
maximal antichain of $\bbB_{\alpha +1}$ included in $\{a \in 
\bbB_{\mathbf k,\alpha +1}:\pi_{\mathbf k,\alpha +1,\alpha}(a) =
0_{\bbB_{\bfk,\alpha}}\}$.

\noindent
2) We say $\mathbf k$ is reasonable \when \, it is reasonable in
   $\alpha$ whenever $\alpha +1 \le \alpha_{\mathbf k}$.

\noindent
3) Let
\mn
\begin{enumerate}
\item[$\bullet$]  $A^1_{\mathbf k,\alpha} = \{f \in \bbN^{\bbB[\mathbf k]}:
f \in \bbN^{\bbB[\mathbf k,\alpha]}$ and 
if $\beta < \alpha$ then $f/D_{\mathbf k} \notin
 \bbN^{\bbB[\mathbf k,\beta]}/D_{\mathbf k}\}$
\sn
\item[$\bullet$]  $A^2_{\mathbf k,\alpha} = \{f/D_{\mathbf k}:f \in  
A^1_{\mathbf k,\alpha}\}$
\sn
\item[$\bullet$]  $A^1_{\mathbf k,< \alpha} = \bigcup\limits_{\beta <
  \alpha} A^1_{\mathbf k,\beta}$ and $A^2_{\mathbf k,< \alpha} =
  \bigcup\limits_{\beta < \alpha} A^2_{\mathbf k,\beta}$, etc.
\end{enumerate}
\mn
4) We say $f$ is reasonable\footnote{Note that ``some $f$ is
  reasonable in $(\bfk,\alpha)$" is close to but not equivalent to
  ``$\bfk$ is reasonable in $\bfh$".}
 in $(\mathbf k,\alpha,\beta)$ when $\alpha <
\beta < \alpha_{\mathbf k}$ and $f \in \bbN^{\bbB[\mathbf k,\beta]}$ and
for some $f' \in \bbN^{\bbB[\mathbf k,\beta]}$, we have $f'/D_{\mathbf k}
= f/D_{\mathbf k}$ and $f'$ is represented
by $\langle a_n:n < \omega \rangle$ and $\pi_{\mathbf k,\alpha}(a_n)=0$
for every $n$ large enough and $a_n \notin D_{\bfk}$ for every $n$.  
If $\beta = \alpha +1$ we may omit it.

\noindent
5) We say $f$ is reasonable in $(\mathbf k,< \alpha)$ when it is
reasonable in $(\mathbf k,\beta,\gamma)$ for some $\beta +1 = \gamma < \alpha$.
\end{definition}

\begin{observation}
\label{m0}
If $\beta < \alpha_{\mathbf k},\mathbf k \in K^{\com}$ and $f \in
\bbN^{\bbB[\mathbf k]}$ is represented by $\langle a_n:n <
\omega\rangle$, \then \, $f \in A^1_{\mathbf k,\le \beta}$ \Iff \,
$\bigcup\limits_{n}(a_n \triangle \pi_{\mathbf k,\beta}(a_n)) \notin D_{\mathbf k}$.
\end{observation}

\begin{PROOF}{\ref{m0}}
By the proof of \ref{x33}.
\end{PROOF}

\begin{claim}
\label{f26}
1) If $\mathbf k \in K^{\cc}_{\alpha +1}$ \then \, there is $\mathbf m \in
K^{\cc}_{\alpha +1}$ such that $\mathbf k \le^{\at}_K \mathbf m,\bbB_{\mathbf m}$ is
complete and $\mathbf m$ is reasonable in $\alpha$ and $\|\bbB_{\mathbf
  m}\| = \|\bbB_{\mathbf k}\|^{\aleph_0}$.

\noindent
2) If $\mathbf k \in K^{\cc}_{\alpha +1}$ 
is reasonable in $\alpha$ and $\mathbf k \le^{\at}_K \mathbf m$ or 
$\mathbf k \le^{\ver}_K \mathbf m$ \then \, $\mathbf m$ is
   reasonable in $\alpha$.

\noindent
3) If $\langle \mathbf k_i:i < \delta\rangle$ is
   $\le^{\ver}_K$-increasing in $K^{\cc}$ and each $\mathbf k_i$ is
   reasonable \then \, there is a $\le^{\ver}_K$-upper bound $\mathbf k$
of cardinality $(\sum\limits_{i} \|\bbB_{\mathbf k_i}\|)^{\aleph_0}$
which is reasonable.

\noindent
4) If $f$ is reasonable in $(\mathbf k,\alpha)$ \then \, it is
   reasonable in $(\mathbf k,< \alpha +1)$.

\noindent
5) If $f \in A^1_{\mathbf k,\alpha}$ \then \, $f$ is reasonable in
   $(\mathbf k,\alpha)$.

\noindent
6) In \ref{f13},(2),(3) if $\mathbf k_i$ is reasonable for every $i <
   \delta$ \then \, so is $\mathbf k$.
\end{claim}

\begin{PROOF}{\ref{f26}}
Straightforward; e.g.:

\noindent
2) Because $\bbB_{\mathbf k} \lessdot \bbB_{\mathbf m,\alpha(\mathbf k)}$,
see Definition \ref{f8}(2) and read Definition \ref{f24}(1).

\noindent
5) Let $\langle a_n:n < \omega\rangle$ represent $f$.

Let $a'_n = \pi_{\mathbf k,\alpha +1} (a_n)$, so $\pi_{\mathbf k,\alpha
  +1,\alpha} (a_n - a'_n) =\pi_{\mathbf k,\alpha +1,\alpha} (a_n) -
  \pi_{\mathbf k,\alpha +1,\alpha}(a'_n) = a'_n - a'_n =0$.
Now $\langle a'_n:n < \omega\rangle$ is an antichain (using
$\pi_{\bfk,\alpha}$) and let $g \in \bbN^{\bbB_\alpha}$ be such that
$g(a'_n)=n$ and $g(1_{\bbB_{\bfk,\alpha}} - \bigcup\limits_{n} a'_n) =
0$.  So $g \in \bbN^{\bbB_{\bfk}}$ and by ``$f \in
A^1_{\bfk,\alpha}$", necessarily $\TV(f \ne g) \in D_{\bfk}$ and
clearly $\TV(f \ne g) \le \bigcup\limits_{n} (a_n-a'_n) \cup \TV(f=0)$
but $\TV(f \ne 0) \in D_{\bfk}$, i.e. $a_0 \notin D_{\bfk}$ by the
assumption on $f$, hence necessarily $b := \bigcup\limits_{n}(a_n -
a'_n) \in D_{\bfk}$.  Now define $f'' \in \bbN^{\bbB_k}$ represented
by $\langle b_n:n < \omega\rangle$ where $b_n = a_n - a'_n$ for $n \ge
1$ and $b_n = 1-\bigcup\limits_{m \ge 1} b_m$ for $n=0$.  So
$\pi_{\bfk,\alpha}(b_n) = 0_{\bbB_{\bfk,n}}$ for $n \ge 1$ and $b_n
\notin D_{\bfk}$ for $n \ge 0$.  Clearly $f''$
is reasonable in $(\mathbf k,\alpha)$.
\end{PROOF}

\begin{claim}
\label{f27}
If (A) then (B) where:
\mn
\begin{enumerate}
\item[(A)]
\begin{enumerate}
\item[(a)]  $\mathbf k \in K^{\cc}_\alpha$
\sn
\item[(b)]  $\beta_n < \beta_{n+1} < \beta =
  \bigcup\limits_{k} \beta_k \le \alpha$
\sn
\item[(c)]  $\mathbf k$ is reasonable in $\beta_n$ for every $n \in \omega$
\end{enumerate}
\sn
\item[(B)]  if $f_1 \in A^1_{\mathbf k,\beta}$, i.e. $f_1 \in 
\bbN^{\bbB[\mathbf k,\beta]},f_1/D_{\mathbf k}
  \notin \cup\{\bbN^{\bbB[\mathbf k,\gamma]}/D_{\mathbf k}:\gamma <
  \beta\}$ \then \, there is $f_2$ such that:
\sn
\begin{enumerate}
\item[(a)]  $f_2 \in \bbN^{\bbB[\mathbf k,\beta]}$
\sn
\item[(b)]  $f_2/D_{\mathbf k} \notin 
\cup\{\bbN^{\bbB[\mathbf k,\gamma]}/D_{\mathbf k}:\gamma < \beta\}$
\sn
\item[(c)]  $f_2/D_{\mathbf k} < f_1/D_{\mathbf k}$
\sn
\item[(d)]   there is $\langle (a_i,k_i):i <
  \omega\rangle$ representing $f_2$ such that:
\sn
\begin{enumerate}
\item[$(\alpha)$]   for each $i$, letting $k(i) = k_i$ we
have $a_i \in \bbB_{\mathbf k,\beta_{k(i)+1}}$ and 
$\pi_{\mathbf k,\beta_{k(i)+1},\beta_{k(i)}} (a_i) =0$
\sn 
\item[$(\beta)$]   for each $\ell$ the set 
$\{i:k_i < \ell\}$ is finite.
\end{enumerate}
\end{enumerate}
\end{enumerate}
\end{claim}

\begin{PROOF}{\ref{f27}}
For each $n$ let $\langle a_{n,\ell}:\ell < \omega\rangle$ be a
maximal antichain of $\bbB_{\beta_n+1}$ such that $\pi_{\mathbf
  k,\beta_n+1,\beta_n}(a_{n,\ell})=0$ for $\ell < \omega$, exists as
$\mathbf k$ is reasonable in $\beta_n$ for every $n \in \omega$, 
see Definition \ref{f24}(2).

Let
\mn
\begin{enumerate}
\item[$(*)_0$]  $(a) \quad \cT_n = \{\eta:\eta \in {}^n \omega$ and
  $\bigcap\limits_{k<n} a_{k,\eta(k)} >0\}$
\sn
\item[${{}}$]  $(b) \quad \cT = \bigcup\limits_{n} \cT_n$.
\end{enumerate}
\mn
Hence
\mn
\begin{enumerate}
\item[$(*)_1$]  $(a) \quad\langle a_\eta:\eta \in \cT_n\rangle$ is a maximal
  antichain of $\bbB_{\beta_{n+1}}$ on which $\pi_{\mathbf k,\beta_n}$
  is zero
\sn
\item[${{}}$]  $(b) \quad \cT$ is a subtree of ${}^{\omega >}\omega$.
\end{enumerate}
\mn
Now choose a sequence $\bar k$ of natural numbers such that:
\mn
\begin{enumerate}
\item[$(*)_2$]  $(a) \quad \bar k = \langle k_s:s \in \cT \rangle$
\sn
\item[${{}}$]  $(b) \quad$ if $\nu \triangleleft \eta$ then $k_\nu <
  k_\eta$
\sn
\item[${{}}$]  $(c) \quad$ if $k_\eta = k_\nu$ then $\eta = \nu$
\sn
\item[$(*)_3$]  let $g_n \in \bbN^{\bbB[\mathbf k,\beta(n)+1]}$ be
  represented by $\langle (a_\eta,k_\eta):\eta \in \cT_n\rangle$, see
  Definition \ref{x20}(4).
\end{enumerate}
\mn
[Why?  Well defined by $(*)_1(a)$.]
\mn
\begin{enumerate}
\item[$(*)_4$]  let $\mathbf M_{\mathbf k} = \{\bar c:\bar c = \langle
  c_\ell:\ell < \omega\rangle$ is a maximal antichain of $\bbB_{\beta_0
  +1}$ disjoint to $D_{\mathbf k}\}$.
\end{enumerate}
\mn
What is the point of $\mathbf M_{\mathbf k}$?  $g_n \in A^1_{\mathbf
  k,\beta(n)}$ hence $\langle g_n/D_{\mathbf k}:n < \omega\rangle$ is
increasing and cofinal in $\cup\{\bbN^{\bbB[\mathbf k,\beta(n)]}/D:n <
\omega\}$ hence if in $\bbN^{\bbB[\mathbf k]}/D_{\mathbf k}$ we have a
definable sequence, the $n$-th try being $g_n/D$, in ``non-standard
places" we have the $g_{\bar c}$'s defined below members of
$A^2_{\mathbf k,\alpha}$ and those are co-initial in it.
\mn
\begin{enumerate}
\item[$(*)_5$]  for each $\bar c \in \mathbf M_{\mathbf k}$ let
\sn
\item[${{}}$]  $(a) \quad S_{\bar c} = \{(\ell,\eta):\ell <
  \omega,\eta \in \cT_\ell$ and $c_\ell \cap a_\eta >0\}$
\sn
\item[${{}}$]  $(b) \quad$ for $(\ell,\eta) \in S_{\bar c}$ let
  $a_{(\ell,\eta)} = c_\ell \cap a_\eta$
\sn
\item[${{}}$]  $(c) \quad g_{\bar c} \in \bbN^{\bbB[\mathbf k,\beta]}$
  be represented by $\langle (a_{(\ell,\eta)},k_\eta):(\ell,\eta) \in
  S_{\bar c}\rangle$
\sn
\item[$(*)_6$]  $g_n$ is $(\mathbf k,\beta_n)$-reasonable.
\end{enumerate}
\mn
[Why?  By the Definition \ref{f24}.]
\mn
\begin{enumerate}
\item[$(*)_7$]  $g_n \in A^1_{\mathbf k,\beta_n +1}$.
\end{enumerate}
\mn
[Why?  Follows from the definition of $g_n$ in $(*)_3$ and the choice
of the $a_\eta$'s and $k_\nu$'s in $(*)_1$ and $(*)_2$.]
\mn
\begin{enumerate}
\item[$(*)_8$]  $g_{\bar c} \in A^1_{\mathbf k,\beta}$ for $\bar c \in
  \mathbf M_{\mathbf k}$.
\end{enumerate}
\mn
[Why?  As $\bar k$ is with no repetition and the definition.]
\mn
\begin{enumerate}
\item[$(*)_9$]  there is $\bar c \in \mathbf M_{\mathbf k}$
such that $g_{\bar c}/D_{\mathbf k} < f_1/D_{\mathbf k}$.
\end{enumerate}
\mn
[Why?  See explanation after $(*)_4$ and \ref{f26}(5).]
\end{PROOF}

\begin{claim}
\label{f39}
If (A) then (B) \where \,:
\mn
\begin{enumerate}
\item[(A)]
\begin{enumerate}
\item[(a)]   $\mathbf k \in K^{\cc}_\alpha$ is $\bbB_{\mathbf k,\alpha}$ infinite
\sn
\item[(b)]  $\beta_n = \beta(n) < \alpha$ is increasing with limit $\alpha$
\sn
\item[(c)]  $\mathbf k$ is reasonable in $\beta_n$
\sn
\item[(d)]  $f \in A^1_{\mathbf k,\alpha}$
\end{enumerate}
\sn
\item[(B)]  there are $\mathbf m,g$ such that:
\sn
\begin{enumerate}
\item[(a)]  $\mathbf k \le^{\at}_K \mathbf m$ and 
$\|\bbB_{\mathbf m}\| = (\|\bbB_{\mathbf k}\|^{\aleph_0})^+$
\sn
\item[(b)]  $g \in A^1_{\mathbf m,\alpha}$
\sn
\item[(c)]  $g/D_{\mathbf m} < f/D_{\mathbf m}$
\sn
\item[(d)]  $g/D_{\mathbf m} \notin \bbN^{\bbB[\mathbf k]}/D_{\mathbf m}$.
\end{enumerate}
\end{enumerate}
\end{claim}

\begin{PROOF}{\ref{f39}}
\Wilog \, $\bbB_{\mathbf k}$ is complete of cardinality $\lambda$.

Let $f_2,\langle (a_n,k_n):n < \omega\rangle$ be as in
\ref{f27} for $f_1 = f$ and let 
$u_n = \{\ell:a_\ell \in \bbB_{\mathbf k,\beta_n}\}$, by \ref{f27}
clearly $u_n$ is finite.

Let $\bbB^0$ be the Boolean algebra extending $\bbB_{\mathbf k}$ generated by
$\bbB_{\mathbf k} \cup \{x_{\varepsilon,n,\ell}:\ell \le n$ and $\varepsilon <
\lambda^+\}$ freely except the equation $x_{\varepsilon,n,\ell} \le
a_n,x_{\varepsilon,n,\ell_1} \cap x_{\varepsilon,n,\ell_2} = 0$ and
$\bigcup\limits_{\ell \le n} x_{\varp,n,\ell} = a_n$ for 
$\varepsilon < \lambda^+,\ell\le n,\ell_1 < \ell_2 \le n$ and let
$\bbB$ be the completion on $\bbB^0$.  Let $g_\varepsilon \in
\bbN^{\bbB}$ be represented by $\langle x_{g_\varepsilon,\ell} :=
\cup\{x_{\varepsilon,n,\ell}:n$ satisfies $\ell \le n\}:\ell <
\omega\rangle$, clearly
\mn
\begin{enumerate}
\item[$(*)_1$]  $g_\varepsilon/D \le f_2/D$ for any $D \in \uf(\bbB)$.
\end{enumerate}
\mn
For $\varepsilon \ne \zeta < \lambda^+$ let $c_{\varepsilon,\zeta} =
\bigcup\limits_{\ell}(x_{g_\varepsilon,\ell} \triangle x_{g_\zeta,\ell})$.

Now
\mn
\begin{enumerate}
\item[$(*)_2$]  $c_{\varepsilon,\zeta} = \bigcup\limits_{n} \,
\bigcup\limits_{\ell \le  n}(x_{\varepsilon,n,\ell} \triangle 
x_{\zeta,n,\ell})$ and $c_{\varp,\zeta} = c_{\zeta,\varp}$.
\end{enumerate}
\mn
[Why?  As $\langle x_{\varepsilon,n,\ell}:\ell \le n\rangle$ is a
partition of $a_n$ and $\langle a_n:n
 < \omega \rangle$ is an antichain of $\bbB$.]

Let $\bbB'$ be the sub-algebra of $\bbB$ generated by $\bbB_{\mathbf k}
\cup Y$ where $Y := \{c_{\varepsilon,\zeta}:\varepsilon < \zeta < \lambda^+\}$
\mn
\begin{enumerate}
\item[$(*)_3$]  we define $\pi^1_n:\bbB_{\mathbf k} \cup Y \rightarrow
  \bbB_{\beta_n}$ by:
\sn
\item[${{}}$]  $\bullet \quad \pi^1_n \rest \bbB_{\mathbf k} =
  \pi_{\mathbf k,\alpha(\mathbf k),\beta(n)}$
\sn
\item[${{}}$]  $\bullet \quad \pi^1_n(c_{\varepsilon,\zeta}) =
  1_{\bbB_{\beta(n)}}$ for $\varepsilon \ne \zeta < \lambda^+$
\sn
\item[$(*)_4$]  $\pi^1_n$ has an extension $\pi^2_n \in
  \Hom(\bbB',\bbB_{\beta(n)})$, necessarily unique
\end{enumerate}
\mn
[Why?  It is enough to show that if $d_0,\dotsc,d_{m-1} \in \bbB_{\mathbf k}$ and
  $\varepsilon_\ell < \zeta_\ell < \lambda^+$ for $\ell < k$ 
and $\sigma(y_0,\dotsc,y_{m-1},x_0,\dotsc,x_{k-1})$ is a Boolean term
and

\[
\bbB \models \sigma(d_0,\dotsc,d_{m-1},c_{\varepsilon_0,\zeta_0},
\dotsc,c_{\varepsilon_{k-1},\zeta_{k-1}})=0
\]

\mn
then $\bbB_{\mathbf k,\beta(n)} \models
\sigma(\pi^1_n(d_0),\dotsc,\pi^1_n(d_{m-1}),
\pi^1_n(c_{\varepsilon_0,\zeta_0},\ldots))=0$.  As
$d_0,\dotsc,d_{m-1} \in \bbB_{\mathbf k}$ and
$\pi^1_n(c_{\varepsilon_\ell,\zeta_\ell}) = 
1_{\bbB_{\mathbf k,\beta(n)}}$ it is sufficient to prove: if $d \in
\bbB_{\mathbf k}$ and $\bbB \models ``d \cap \bigcap\limits_{\ell < k}
c_{\varepsilon_\ell,\zeta_\ell} = 0"$ then 
$\bbB_{\mathbf k,\beta(n)} \models ``(\pi_{\mathbf k,\alpha(\mathbf k),
\beta(k)}(d)) = 0"$.

Now if for some $\ell \ge 1,\ell \notin u_n,d \cap a_\ell > 0$ the assumption
does not hold and otherwise, necessarily $d \le \bigcup\limits_{\ell
  \in u_n} a_\ell$ hence the conclusion holds.  So indeed $\pi^2_n$ is
well defined.]
\mn
\begin{enumerate}
\item[$(*)_5$]  $\pi^2_n = \pi_{\mathbf k,\beta(n+1),\beta(n)} \circ
 \pi^2_{n+1}$.
\end{enumerate}
\mn
[Why?  See the definition of $\pi^1_m$ recalling \ref{f17}(4).]
\mn
\begin{enumerate}
\item[$(*)_6$]  there is $\mathbf n$ such that $\mathbf k \le^{\at}_K
 \mathbf n$ and $\bbB_{\mathbf n} = \bbB',\pi_{\mathbf n,\alpha(\mathbf
 k),\beta(n)} = \pi^2_n$ hence $D_{\mathbf n} \supseteq Y$.
\end{enumerate}
\mn
[Why?  Check the definitions.]
\mn
\begin{enumerate}
\item[$(*)_7$]  there is $\mathbf m$ such that $\mathbf n \le^{\wa}_K
 \mathbf m,\bbB_{\mathbf m} = \bbB$ hence $\mathbf k \le^{\wa}_K \mathbf m$
\end{enumerate}
\mn
[Why?  By claim \ref{f16}(1),(3), the ``hence" by \ref{f16}(5)
recalling that $\bbB_{\bfk} \lessdot \bbB$ by the choice of $\bbB$ and
$\bbB_{\bfk} \subseteq \bbB' \subseteq \bbB$ by the choice of $\bbB'$.]
\mn
\begin{enumerate}
\item[$(*)_8$]  there is $\varepsilon < \lambda^+$ such that
 $g_\varepsilon \in A^1_{\mathbf n,\alpha(\mathbf k)}$.
\end{enumerate}
\mn
[Why?  As $A^1_{\mathbf n,< \alpha(\mathbf k)}$ has cardinality $\le \lambda$.]
\mn
\begin{enumerate}
\item[$(*)_9$]   $\mathbf m$ is as required.
\end{enumerate}
\mn
[By $(*)_7 + (*)_8$ and \ref{f16}(2).]
\end{PROOF}

\begin{observation}
\label{f40}
In claim \ref{f39} we can demand $\|\bbB_{\bfm}\| = 
\|\bbB_{\mathbf k}\|^{\aleph_0}$.
\end{observation}

\begin{PROOF}{\ref{f40}}
By the L\"owenheim-Skolem-Tarski argument.
\end{PROOF}

\begin{claim}
\label{f42}
Assume $\mathbf k \in K^{\cc}_\alpha$ and $\cf(\alpha) > \aleph_0$.

If $f \in \bbN^{\bbB[\mathbf k]}$ and $f/D \notin 
\bbN^{\bbB[\mathbf k,\beta]}/D_{\mathbf k}$ for $\beta < \alpha$, 
\then \, for some $\mathbf m$ and $g$ we have:
\mn
\begin{enumerate}
\item[$(*)$]  $(a) \quad \mathbf k \le^{\at}_K \mathbf m$
\sn
\item[${{}}$]  $(b) \quad g \in \bbN^{\bbB[\mathbf m]}$
\sn
\item[${{}}$]  $(c) \quad g/D_{\mathbf m} < f/D_{\mathbf m}$
\sn
\item[${{}}$]  $(d) \quad h/D_{\mathbf m} < g/D_{\mathbf m}$ \when \, $h
  \in \bbN^{\bbB[\mathbf k,\beta]}$ for some $\beta < \alpha$
\sn
\item[${{}}$]  $(e) \quad |\bbB_{\mathbf m}| \le |\bbB_{\mathbf k}|$.
\end{enumerate}
\end{claim}

\begin{PROOF}{\ref{f42}}
Like the proof of \ref{f39} only simpler and shorter.  Let $\bar a =
\langle a_n:n < \omega\rangle$ represent $f,\lambda = 
\|\bbB_{\mathbf k}\|$.  By \ref{f26}(4), \wilog \, 
$f$ is reasonable in $(\bfk,\alpha)$; let 
$\{x_{\varepsilon,n,\ell}:\varepsilon < \lambda^+,\ell \le n\},
\bbB^0,\bbB,Y,\bbB'$ be as there and define $\pi_1:\bbB_{\mathbf k} 
\cup Y \rightarrow \bbB_{\mathbf k,\alpha}$ as there.  $\pi_1 \rest 
\bbB_{\mathbf k,\alpha} = \pi_{\mathbf k,\alpha +1,\alpha},
\pi_1(x_{\varepsilon,\zeta}) = 1_{\bbB_{\mathbf k,\alpha}}$ 
for $\varepsilon < \zeta < \lambda^+$.

The rest is as there. 
\end{PROOF}

\begin{claim}
\label{f34}
If $\mathbf k \in K^{\com}_\alpha,\lambda \ge \|\bbB_{\mathbf k}\| +
2^{\aleph_0}$ and $p(x)$ is a type in the model
  $\bbN^{\bbB[\mathbf k]}/D_{\mathbf k}$ \then \, for some $\mathbf m \in
  K_{\alpha +1}$ we have $\mathbf k \le^{\ver}_K \mathbf m$ and $p(x)$ is
  realized in $\bbN^{\bbB[\mathbf m]}/D_{\mathbf m}$.
\end{claim}

\begin{PROOF}{\ref{f34}}
Easy.
\end{PROOF}

\noindent
Having established all these statements, we can prove now the main
result of this paper:
\begin{theorem}  
\label{k6}
For any infinite cardinal $\lambda$, for 
some regular ultrafilter $D$ on $\lambda$ we have $\upf(D) = \cC$
\Iff \,:
\mn
\begin{enumerate}
\item[$(*)$] 
\begin{enumerate}
\item[(a)]  $\cC$ is a set of cardinals $\le 2^\lambda$
\sn
\item[(b)]   $\mu = \mu^{\aleph_0}$ whenever $\mu \in \cC$
\sn
\item[(c)]  $2^\lambda \in \cC$ hence is the maximal member of $\cC$.
\end{enumerate}
\end{enumerate}
\end{theorem}

\begin{PROOF}{\ref{k6}}
The implication $\Rightarrow$, we already know, so we shall deal with the
$\Leftarrow$ implication; the proof relies on earlier definitions and
claims so the reader can return to a second reading.

Let $\langle \lambda_\alpha:\alpha \le \alpha(*)\rangle$ list $\cC$ in
increasing order.  Let $S = \{\alpha:\alpha \le \alpha(*)+1$ and 
$\cf(\alpha) \ne \aleph_0\}$.
We choose $\mathbf k_\alpha$ by induction on $\alpha \in S \cap
(\alpha(*)+2)$ such that:
\mn
\begin{enumerate}
\item[$\boxplus$]
\begin{enumerate}
\item[(a)]  $\mathbf k_\alpha \in K^{\com}_\alpha \cap
  K^{\cc}_\alpha$, see Definition \ref{f2}, \ref{f4}(1A)
\sn
\item[(b)]  $\mathbf k_\beta \le^{\ver}_K \mathbf k_\alpha$
 for $\beta \in \alpha \cap S$, see \ref{f8}(2)(B)
\sn
\item[(c)]  if $f \in A^1_{\bfk,\beta},\beta < \alpha$ then
$\lambda_\beta = |\{g/D_{\bfk}:g \in \bbN^{\bbB},g/D < f(D)\}|$
\sn
\item[(d)]  if $\cf(\alpha) > \aleph_0$ then
$\bbB_{\mathbf k_\alpha} = \cup\{\bbB_{\mathbf k_\beta}:\beta < \alpha\}$
\sn
\item[(e)]  $\mathbf k_\alpha$ is reasonable (see Definition \ref{f24})
\sn
\item[(f)]  the set $A^1_{\mathbf k,\beta}$ has cardinality
  $\lambda_\beta$.
\end{enumerate}
\end{enumerate}
\medskip

\noindent
\underline{Case 1}:  For $\alpha = 0$

$\bbB_{\mathbf k_0}$ is the trivial Boolean Algebra, so really there is
nothing to prove.
\medskip

\noindent
\underline{Case 2}:   $\cf(\alpha) > \aleph_0$

Use \ref{f13}(2),(3) to find $\mathbf k_\alpha$ satisfying clauses
(a),(b),(c),(d).  Now $\mathbf k_\alpha$ satisfies clause (e) by \ref{f26}(6).
\medskip

\noindent
\underline{Case 3}:  $\alpha = \beta +1$

We choose $\mathbf k_{\beta,i}$ by induction for $i \le \lambda_\beta$
such that
\mn
\begin{enumerate}
\item[$(*)$]  
\begin{enumerate}
\item[(a)]
\begin{enumerate}
\item[$(\alpha)$]  if $\beta \in S$ then 
$\mathbf k_\beta \le^{\ver}_K \mathbf k_{\beta,i} \in K^{\com}_\alpha \cap
  K^{\cc}_\alpha$
\sn
\item[$(\beta)$]  if $\beta \notin S$ then
$\gamma \in \beta \cap S \Rightarrow \mathbf k_\gamma \le^{\ver}_K
\mathbf k_{\beta,i} \in K^{\com}_\alpha \cap K^{\cc}_\alpha$
\sn
\item[$(\gamma)$]  if $i=0$ then there is
  $g \in \bbN^{\bbB[\mathbf k_{\beta,i}]}$ such that
$g/D_{\mathbf k_{\beta,i}} \notin \{f/D_{\mathbf k_{\beta,i}}:
f \in \bbN^{\bbB[\mathbf k_\gamma]}$ for some $\gamma \in \alpha \cap S\}$
\sn
\item[$(\delta)$]  $\bbB_{\mathbf k_{\beta,i}}$ is infinite
\end{enumerate}
\sn
\item[(b)]  $\langle \mathbf k_{\beta,j}:j \le i\rangle$ is
$\le^{\at}_K$-increasing continuous
\sn
\item[(c)]  $\bbB_{\mathbf k_{\beta,i}}$ has cardinality $\le \lambda_\beta$
\sn
\item[(d)]   if $i=j+1$:
\sn
\begin{enumerate}
\item[$(\alpha)$]   bookkeeping gives us $g_{\beta,j} \in
\bbN^{\bbB[\mathbf k_{\beta,i}]}$ such that 
$g_{\beta,j}/D_{\mathbf k_{\beta,j}} \notin 
\cup\{\bbN^{\bbB[\gamma,\mathbf k_\beta]}:\gamma \in \alpha \cap S\}$
\sn
\item[$(\beta)$]   there is $f_{\beta,j} \in
\bbN^{\bbB[\mathbf k_{\beta,i}]}$ such that $f_{\beta,j}/D_{\mathbf
k_{\beta,i}} < g_{\beta,j}/D_{\mathbf k_{\beta,i}}$ and
$f_{\beta,j}/D_{\mathbf k_{\beta,i}} \notin
\cup\{\bbN^{\bbB[\gamma,\mathbf k_\beta]}:\gamma \in \alpha \cap S\}$
\end{enumerate}
\sn
\item[(e)]  if $i < \lambda_\beta$ and $g$ satisfies 
$(d)(\alpha)$ then for some $i_1 \in [i,\lambda_\beta],g_{\beta,i_1} =g$
\sn
\item[(f)]  if $i=j+1$ then $\bbB_{\mathbf k_{\beta,i}}$
  is complete and reasonable.
\end{enumerate}
\end{enumerate}
\mn
Note that by \ref{f26}(1) we can take care of clause (f), so we shall
ignore it.

For $i = 0$ we use \ref{f26}(1) if $\beta \in S$ and we use 
\ref{f13} if $\beta \notin S$.

For $i$ limit use \ref{f13}(1).

For $i=j+1,\cf(j) > \aleph_0$ use the claim \ref{f42}.

If $i=j+1,\cf(j) = \aleph_0$ use the claim \ref{f39}.

For $i=j+1,\cf(j)=1$ we use \ref{f26}(1).

Having carried the induction on $i \le \lambda_\beta$ let
$\mathbf k_\alpha = \mathbf k_{\beta,\lambda_\beta}$.  In particular
$\bbB_{\mathbf k_\beta,\lambda_\beta}$ is complete as $\lambda_\beta =
\sup\{i < \lambda_\beta:\bbB_{\mathbf k,i}$ is complete$\}$ by clause
(f) and $\cf(\lambda_\beta) > \aleph_0$ as $\lambda_\beta =
\lambda^{\aleph_0}_\beta$. 

Having carried the induction on $\alpha \le \alpha(*) +1$ clearly the pair
$(\bbB_{\mathbf k_{\alpha(*)+1}},D_{\mathbf k_{\alpha(*)+1}})$ is almost
as required.  That is, (see \cite[Ch.VI,\S3]{Sh:c}) we know that for
some regular filter $D_*$ on $\cP(I)$, there is a homomorphism $\mathbf j$ from
the Boolean Algebra $\cP(I)$ onto $\bbB_{\mathbf k_{\alpha(*)+1}}$ and
let $D = \{A \subseteq \lambda:\mathbf j(A) \in D_{\mathbf k_{\alpha(*)+1}}\}$.
\end{PROOF}
\newpage


\bibliographystyle{amsalpha}
\bibliography{shlhetal}

\end{document}